\documentclass{article}

\usepackage{amsmath, amssymb, amsfonts, amsthm}
\usepackage{mathtools}
\usepackage{geometry}
\usepackage{graphicx}
\usepackage{bm}
\usepackage{caption}
\usepackage{subcaption}
\usepackage{booktabs}
\usepackage{multirow}
\usepackage{array}
\usepackage{algorithm}
\usepackage{algorithmicx}
\usepackage{algpseudocode}
\usepackage{hyperref}
\usepackage{cleveref}
\usepackage{placeins}

\DeclareMathOperator{\tr}{tr}
\DeclareMathOperator{\diag}{diag}

\newcommand{\R}{\mathbb{R}}
\newcommand{\E}{\mathbb{E}}
\newcommand{\norm}[1]{\left\|#1\right\|}

\newcommand{\dd}{\mathrm{d}}
\newcommand{\bs}[1]{\boldsymbol{#1}}
\newcommand{\x}{\mathbf{x}}
\renewcommand{\u}{\mathbf{u}}
\newcommand{\z}{\mathbf{z}}
\newcommand{\bX}{\mathbf{X}}
\newcommand{\KL}[2]{D_{\mathrm{KL}}\!\left(#1\,\|\,#2\right)}

\newtheorem{proposition}{Proposition}

\title{Scalable Multi-Agent Gas Source Localization via\\
       Continuous-Discrete Stochastic Optimal Control}

\author{
  Achim Basermann\textsuperscript{1}\quad
  André Carlon\textsuperscript{2,4}\thanks{Correspondence to André Carlon at
    \href{mailto:carlon@uq.rwth-aachen.de}{\nolinkurl{carlon@uq.rwth-aachen.de}} or
    \href{mailto:agcarlon@gmail.com}{\nolinkurl{agcarlon@gmail.com}}.}\quad
  Philipp Knechtges\textsuperscript{1}\\[0.2em]
  Raúl Tempone\textsuperscript{3}\quad
  Felix Terhag\textsuperscript{1,2}\\[0.6em]
  \parbox{0.95\textwidth}{\centering\small
    \textsuperscript{1}Institute of Software Technology, High-Performance Computing,
    German Aerospace Center (DLR), Cologne, Germany\endgraf
    \textsuperscript{2}Chair of Mathematics for Uncertainty Quantification,
    RWTH Aachen University, Aachen, Germany
    \textsuperscript{3}Computer, Electrical and Mathematical Sciences and Engineering Division (CEMSE),
    King Abdullah University of Science and Technology (KAUST), Thuwal, Saudi Arabia\endgraf
    \textsuperscript{4}Center for Optimization and Reliability Engineering,
    UFSC, Florianópolis, Brazil\endgraf
  }
}

\date{}

\begin{document}
\maketitle

\begin{abstract}
Localizing a gas source with multiple unmanned aerial vehicles (UAVs) requires coordinating informative measurements while controlling motion cost.
We present a stochastic control framework that combines continuous-time feedback with discrete Bayesian updates of source location and log-intensity.
A Hamilton--Jacobi--Bellman (HJB) equation balances wind-relative control effort, crosswind patrolling, and collision and no-fly-zone penalties between measurements.
At measurement times, a linearized-Gaussian approximation supplies an expected information gain reward and a predicted covariance update.
A Bellman fixed point over one measurement interval accounts for repeated cycles under this reduced belief model, and the policy is recomputed after each actual observation.
To limit computational cost, we combine one-body or pairwise HJB solutions and use wind-aligned covariance coordinates, yielding five- or seven-dimensional subproblems independent of team size.
In repeated simulation trials, one-body HJB control achieves success rates within four percentage points of a greedy information-gain baseline while using less control effort.
In a challenging three-UAV scenario, pairwise planning increases the fraction of trials localized within $5\,\mathrm{m}$ from $71\%$ to $93\%$, at approximately $16$ times the median runtime.
Weaker discounting improves accuracy at increased motion cost, and experiments with two to four UAVs demonstrate reuse of the fixed-dimensional pairwise planner.
These results establish practical accuracy--effort--computation trade-offs within the tested steady-plume model.
\end{abstract}


\section{Introduction}
\label{sec:introduction}
Gas source localization (GSL) is a fundamental task in environmental monitoring and industrial safety, including pollutant monitoring and leak detection~\cite{francis2022gslReview}.
UAVs provide a natural platform for this task: they can rapidly cover large areas, adapt trajectories online, and reach otherwise inaccessible or dangerous environments.

Despite significant progress, GSL remains a challenging estimation and control problem.
Measurements are sparse, noisy, and indirect, and wind-driven advection introduces strong spatial structure.
When source intensity is unknown, different combinations of location and intensity can produce similar concentrations at a measurement position.
The choice of measurement locations therefore affects how well the two can be distinguished.
The problem becomes substantially more complex in multi-agent scenarios.
While multiple UAVs offer increased sensing coverage, the controller must choose each drone's motion while accounting for the information collected by the other drones.
The state dimension grows linearly with the number of agents and makes grid-based HJB planning exponentially more expensive as the team size increases.

Information-driven search addresses the sensing problem by maintaining a posterior over source parameters and steering toward informative measurement locations~\cite{vergassola2007infotaxis,hutchinson2019control,chen2021dcee}.
Continuous trajectory optimization and informative path planning also support source localization and related sensing tasks~\cite{charrow2015information,jin2024sense,jin2025cumulative}.
The difficulty addressed here is to combine continuous feedback, repeated Bayesian updates, joint inference of source location and intensity, and computationally manageable coordination of several UAVs.

Our approach separates continuous motion from discrete learning.
Between measurements, a Hamilton--Jacobi--Bellman (HJB) equation balances motion and operational costs against the information expected at the next observation.
At measurement times, an information reward and a predicted covariance update connect consecutive planning intervals.
A Bellman fixed point over one interval represents repeated cycles without prescribing a finite terminal horizon.
This architecture follows the continuous-time control and discrete-observation viewpoint of Bayer et al.~\cite{bayer2024continuous}.

The central computational choice is to use a reduced belief model for planning while updating source location and log-intensity jointly from actual measurements through a Laplace approximation~\cite{laplace_approx,long2013fast}.
The planner holds the current parameter estimates fixed and predicts future covariance changes using a block-diagonal Fisher approximation.
One-body or pairwise HJB subproblems then provide feedbacks that can be reused across the team, and a wind-aligned covariance parameterization removes a further grid dimension.
The resulting controller is an approximation to joint belief-state control; its purpose is to retain useful information-driven feedback at a manageable computational cost.

The main contributions of this work are:
\begin{enumerate}
  \item A continuous-time stochastic control formulation for multi-drone GSL with discrete information reward and Bayesian update while inferring source location and log-intensity jointly.
  \item A Bellman fixed point that approximates repeated replanning beyond the next measurement without imposing a finite terminal horizon.
  \item Fixed-dimensional one-body and pairwise feedback surrogates, together with a wind-aligned belief parameterization, that reduce the dimensionality of online HJB replanning.
  \item A simulation study quantifying the accuracy, control-effort, and runtime trade-offs of these approximations, including comparisons with patrol and greedy EIG control and sensitivity to discount rate and team size.
\end{enumerate}

\section{Related Work}
\label{sec:related_work}

Many strategies have been proposed to address the challenges of gas source localization, ranging from biologically inspired heuristics to formal stochastic control frameworks \cite{francis2022gslReview}.

\subsection{Bio-Inspired and Heuristic Search}

Research on GSL has long been inspired by biological observations.
Flying insects navigating toward wind-borne odors often proceed upwind using counterturning patterns: when signal is absent, they perform extended ``casting'' or shorter ``zigzagging'' crosswind maneuvers \cite{kennedy83zigzagging}.
Many robotic algorithms directly adopt these strategies, for example crosswind zigzagging \cite{li2001crosswindimportance} or chemotaxis \cite{russel03chemotaxis}.
Purely reactive versions are attractive for their simplicity and scalability, but do not explicitly represent source uncertainty.

\subsection{Information-Theoretic Planning}

A more principled class of approaches maintains a Bayesian belief over the source and selects actions to maximize expected information gain (EIG), a classical criterion in Bayesian optimal experimental design \cite{lindley1956measurement,chaloner1995bayesian}.
A prominent example is Infotaxis \cite{vergassola2007infotaxis}, which combines a Bayesian posterior with an entropy-reduction reward while retaining the biological intuition of zigzagging and casting.
Many follow-ups use related criteria such as Kullback--Leibler divergence between prior and posterior \cite{hutchinson2019control}, or cast the problem as a partially observed Markov decision process (POMDP) with information rewards \cite{rahbar2019odorsource}.
Several prominent schemes evaluate a discrete set of candidate moves in a greedy one-step lookahead fashion, which can miss longer-term strategies and may stall when the immediate information reward is zero.
Chen et al.\ \cite{chen2021dcee} introduce a Dual Control for Exploitation and Exploration (DCEE) approach that explicitly balances moving toward the current source estimate against reducing uncertainty; while more principled, it still discretizes the action space and uses a myopic horizon.
Continuous alternatives include trajectory refinement for active information gathering \cite{charrow2015information}, sense-in-motion source seeking \cite{jin2024sense}, and cumulative informative path planning for gas-source localization \cite{jin2025cumulative}.
The continuous-time stochastic control framework with discrete-time partial observations of Bayer et al.~\cite{bayer2024continuous} is closely related in control architecture: it formulates partially observed control through belief-state dynamics, Bayesian updates at observation times, and HJB equations between observations.
Our work uses the same continuous-time/discrete-observation viewpoint, but specializes it to gas-source localization with a differentiable plume model, a Fisher-information approximation to EIG, joint source-location/intensity inference, and low-order multi-UAV planning surrogates.

\subsection{Learning-Based Methods}

Deep reinforcement learning has also been applied to GSL.
Zhao et al.\ \cite{zhao2022rl} train a DQN-based policy for single-robot search in a turbulent plume.
In adjacent source-localization settings, Le et al.\ \cite{le25multi_agent_rl} propose a multi-agent RL framework with value decomposition networks to encourage cooperation, albeit in a different sensing modality from gas-plume localization.
Data-driven policies can handle complex environmental physics learned from experience, but they require extensive training data and do not provide problem-specific optimality or safety guarantees in the setting considered here.

\subsection{Multi-Robot Coordination}

In the multi-robot setting, a wide range of coordination strategies have been explored.
A first distinction is between centralized strategies, which fuse all robots' measurements into a common posterior, and decentralized ones, where each agent maintains its own belief and shares information intermittently.
Heuristic approaches often minimize team travel distance or balance spatial coverage while sharing measurements, without explicitly maximizing joint information gain \cite{rahbar2020multirobot}.
Formation-based methods maintain a fixed agent geometry \cite{ristic2017formation}, while swarm-inspired coordination draws on bio-inspired behaviors such as particle swarm optimization \cite{bourne2019bioinspired} or ant-colony optimization \cite{meng2006aco}.

On the information-theoretic side, more recent approaches exploit mutual information between sensor measurements and source estimates to coordinate multiple robots \cite{shaffer2025multirobot}, but often rely on distributed or approximate coordination rather than solving a fully joint continuous-time planning problem over the coupled team state.
Bourne et al.~\cite{bourne2020outdoor} develop decentralized, trajectory-based information-theoretic control for multi-agent gas-leak localization, while Wiedemann et al.~\cite{wiedemann2019coop} combine probabilistic, model-based localization with cooperative experiments and explicit model uncertainty.
Related distributed source-seeking methods also provide rigorous coordination analyses in adjacent sensing settings \cite{zhang2023distributed}.
Our distinction from these methods lies primarily in the control formulation; it uses a centralized belief and approximate low-order HJB feedbacks.

The proposed framework combines centralized joint estimation with continuous HJB feedback and a repeated-cycle information reward.
Its computational reduction comes from solving one-body or two-body control problems and reusing their feedbacks across the team.
The comparisons below therefore assess the complete approximate controller, including the cost of retaining pairwise interactions.

\section{Problem Formulation}
\label{sec:problem}

We now specify the plume, sensing, belief, and motion models used by the proposed control framework.
We consider the problem of localizing an unknown pollutant source using a team of $M$ drones whose nominal operating region is a bounded 2-D domain $\Omega\subset\R^2$.
The source is fixed at an unknown location $\bs\theta\in\Omega$ and emits pollutant at rate $q>0$.
A steady wind $\bs w = w_s[\cos\alpha,\,\sin\alpha]^\top$ advects the plume with speed $w_s$ and direction $\alpha$, where $\alpha$ is measured counterclockwise from the positive $x_1$-axis.
Thus, $\bs w$ points downwind and $-\bs w$ points upwind.

\subsection{Concentration Model}
\label{sec:plume_model}
Rather than using a full atmospheric-dispersion model, we use an analytically tractable sliced screened plume surrogate derived from an advection-diffusion-decay equation with a point source.
Let $d_s>0$ denote the pollutant diffusivity, let $\tau_s>0$ denote its decay lifetime, and write the unknown source rate as
$q=q_{\mathrm{ref}}e^\zeta$, where $q_{\mathrm{ref}}>0$ is a fixed reference rate and $\zeta=\log(q/q_{\mathrm{ref}})$ is dimensionless.
For a horizontal position $\x\in\R^2$, define its three-dimensional embedding $\widetilde{\x}:=(\x^\top,0)^\top$, and define $\widetilde{\bs\theta}$ and $\widetilde{\bs w}$ analogously.
Starting from the steady three-dimensional equation
\begin{equation}
-d_s \Delta C + \widetilde{\bs w}\cdot \nabla C + \tau_s^{-1} C
= q_{\mathrm{ref}}e^\zeta \delta_3(\widetilde{\x}-\widetilde{\bs\theta})
\qquad \text{in } \R^3,
\end{equation}
% REVIEW C3-plume BEGIN
where $\delta_3$ denotes the three-dimensional Dirac distribution.
Writing $\bs y=\widetilde{\x}-\widetilde{\bs\theta}$ and $C=e^{\widetilde{\bs w}\cdot\bs y/(2d_s)}U$ reduces the operator to $-d_s\Delta U+(\tau_s^{-1}+\norm{\bs w}^2/(4d_s))U$.
Its decaying point-source solution is $U(\bs y)=q\exp(-\kappa\norm{\bs y})/(4\pi d_s\norm{\bs y})$; see \cite{Stockie2011} for general dispersion-modeling background.
Restricting to the horizontal plane and choosing a core regularization of the radial terms gives the surrogate
% REVIEW C3-plume END
\begin{equation}
\label{eq:plume_model}
C(\x,\bs\theta,\zeta)=
\frac{q_{\mathrm{ref}}e^{\zeta}}{4\pi d_s r_c(\x,\bs\theta)}
\exp\!\left(-\kappa r_c(\x,\bs\theta)\right)
\exp\!\left(\frac{\bs w\cdot(\x-\bs\theta)}{2d_s}\right),
\end{equation}
where $\varepsilon_c>0$ is a core-regularization length and
\begin{equation}
r_c(\x,\bs\theta)=\sqrt{\norm{\x-\bs\theta}^2+\varepsilon_c^2},
\qquad
\kappa^2=\frac{1}{d_s\tau_s}+\frac{\norm{\bs w}^2}{4d_s^2}.
\end{equation}
Thus, $r_c$ is the regularized source distance and $\kappa>0$ is the inverse screening length.
This steady surrogate captures the wind-aligned anisotropy used by the planner while neglecting fine-scale turbulent fluctuations, and it provides closed-form concentration evaluation, analytic sensitivities with respect to $\bs\theta$ and $\zeta$, and inexpensive vectorized evaluation for repeated Fisher-information and HJB computations.
Analogous closed-form steady plume surrogates could be substituted without changing the estimator or planner pipeline.

The experiments express the source rate $q$ in grams per second and use $q_{\mathrm{ref}}=1\,\mathrm{g\,s^{-1}}$.
Because $q$ is difficult to measure independently, we infer the log-intensity $\zeta$ jointly with the source location (see \Cref{sec:belief_update}).
At planning time $t$, the current intensity estimate is $\hat q_t=q_{\mathrm{ref}}e^{\hat\zeta_t}$, so $C(\x,\bs\theta,\hat\zeta_t)=e^{\hat\zeta_t} C(\x,\bs\theta,0)$.

\subsection{Sensing and Dynamics}
\label{sec:sensing_dynamics}
% REVIEW C1-initial-measurement BEGIN
Measurements are collected at times $t_k=t_0+k\Delta$, $k=0,1,\ldots$, with interval $\Delta>0$.
The reported experiments include the observation at $t_0$.
% REVIEW C1-initial-measurement END
At each $t_k$ drone $i$ records
\begin{equation}
\label{eq:measurement_model}
y_{ik}=C\!\left(\x_i(t_k),\bs\theta,\zeta\right)+\epsilon_{ik},\qquad
\epsilon_{ik}\sim\mathcal N(0,\sigma_e^2),
\end{equation}
independently across drones and time steps, where both the source location $\bs\theta$ and the log-intensity $\zeta$ are unknown.
% REVIEW C1-dataset BEGIN
The cumulative dataset after this observation is $D_k=\{(\x_i(t_j),y_{ij}) : i=1,\ldots,M,\; j=0,\ldots,k\}$, with $D_{-1}=\emptyset$ before the initial observation.
% REVIEW C1-dataset END

Between measurements, each drone evolves in $\R^2$ as
\begin{equation}
\label{eq:dynamics}
\dd\x_i = \u_i\,\dd t + \sigma_d\,\dd W_i,\qquad
\norm{\u_i}\le v_{\max},
\end{equation}
where $\sigma_d\ge 0$ is a diffusion coefficient and $W_i$ are independent standard Brownian motions.
We assume a filtered probability space carrying the mutually independent Brownian motions and measurement errors, with drone positions observed continuously and without error.
Let $\mathcal A$ denote the set of admissible controls: controls that are progressively measurable with respect to the filtration generated by the continuously observed positions and measurements received up to the current time, and that satisfy $\norm{\u_i(t)}\le v_{\max}$ almost surely.
The source parameters are static, $\dd\bs\theta=0$ and $\dd\zeta=0$, and are random only through the Bayesian prior and simulation ensemble.
Here $\u_i$ is commanded ground velocity and $v_{\max}$ is a ground-speed bound; consequently, $\norm{\u_i-\bs w}^2$ below is a heuristic air-relative control-effort proxy rather than a complete actuator model.
The bounded set $\Omega$ denotes the nominal operating region.
Trajectories are initialized in $\Omega$.
We do not impose reflecting or absorbing boundary conditions in the physical SDE.
The HJB is solved on a bounded computational box with projection of off-grid semi-Lagrangian foot points; in simulation, physical positions are projected back into $\Omega$ only as a safeguard.
The frequency and duration of physical projection are treated as numerical diagnostics in the experiments.

The diffusion coefficient $\sigma_d$ serves a dual purpose.
At face value it models physical perturbations such as wind gusts or actuator noise.
It also introduces a diffusion term in the HJB equation that smooths the value function over position.
This smoothing may make the policy less sensitive to sharp gradients of the predicted concentration field and can therefore be interpreted as a form of regularization against model mismatch.
We do not, however, claim or test a quantitative robustness guarantee with respect to $\sigma_d$.

\subsection{Gaussian Belief and Control Objective}
\label{sec:belief_objective}
Define the stacked team position and control by
\begin{equation}
\bX_t=[\x_1(t)^\top,\dots,\x_M(t)^\top]^\top\in\R^{2M},
\qquad
\u_t=[\u_1(t)^\top,\dots,\u_M(t)^\top]^\top\in\R^{2M}.
\end{equation}
The centralized estimator maintains a Gaussian approximation to the posterior over the full parameter vector $\bs\xi=(\bs\theta,\zeta)\in\R^3$.
The location and log-intensity components $(\bs\mu_t,\hat\zeta_t)$ of the constrained joint maximum a posteriori (MAP) estimate are the means of this Gaussian approximation.
Its location covariance and log-intensity variance are $\bs\Sigma_t$ and $\sigma_{\zeta,t}^2$, respectively.
The estimator retains location--intensity dependence when computing these marginals; the planner uses a further factorization described in \Cref{sec:eig_computation}.
These marginal quantities are refreshed at each $t_k$ and held fixed between updates.
We seek a feedback control policy $\u^\star(\bX,t)$ that drives the UAVs to localize $\bs\theta$ quickly while balancing motion cost, patrolling, collision avoidance, and zone constraints.

Conceptually, the mission-level localization objective is to get $\bs\mu_t$ as close as possible to the true source $\bs\theta$.
The belief covariance trace $\tr(\bs\Sigma_t)$ measures the local concentration of the Gaussian approximation around its MAP center.
A small trace indicates high local posterior confidence, but does not by itself imply that the MAP estimate is close to the true source.
We therefore combine the trace with measurement-fit and intensity-uncertainty checks in the stopping rule (\Cref{sec:stopping_criterion}).

We assume $w_s$, $\alpha$, $d_s$, and $\tau_s$ are known and constant (known plume physics), and that all drones share the same $\sigma_e^2$ (homogeneous sensors).
The assumption of constant, spatially uniform wind is an idealization.
The Gaussian belief over $\bs\theta$ should be interpreted as an unconstrained approximation, even though the physical operating region is bounded.
All agents communicate without delay (centralized control).

\subsection{Expected Information Gain as a Control Objective}
\label{sec:eig_objective}

Let $\pi_t^-$ denote the current belief over $\bs\xi$, and let $Y_{k+1}$ denote the measurement batch that would be collected at a candidate team configuration $\bX$.
After observing $Y_{k+1}$, Bayes' rule gives the posterior $\pi_t^+(\cdot\mid Y_{k+1},\bX)$.
The cycle-wise expected information gain associated with $\bX$ is
\begin{equation}
\label{eq:cycle_eig_def}
\operatorname{EIG}(\bX;\pi_t^-)
:=
\E_{Y_{k+1}\mid \bX,\pi_t^-}
\!\left[
\KL{\pi_t^+(\cdot\mid Y_{k+1},\bX)}{\pi_t^-}
\right].
\end{equation}
Equivalently, $\operatorname{EIG}(\bX;\pi_t^-)$ is the mutual information between the unknown source parameters and the prospective measurement batch, conditioned on the current belief and candidate configuration.
It therefore rewards configurations whose measurements are expected to distinguish between parameter values that remain plausible under the current belief.

We use EIG as a discrete reward at measurement times.
This connects estimation and control: the UAVs are encouraged to reach informative configurations while balancing the motion and operational costs introduced below.
For a fixed $\bX$, EIG is a deterministic design criterion because the expectation in \eqref{eq:cycle_eig_def} averages over the prospective measurement outcomes.
Exact evaluation is nevertheless generally intractable for the nonlinear plume model.
\Cref{sec:eig_computation} introduces the linearized-Gaussian approximation used to obtain a closed-form planner reward and predicted covariance update.

\section{Method}
\label{sec:method}

The proposed controller treats the interval between two consecutive measurement times as one cycle of a continuous-discrete stochastic optimal control problem.
Within a cycle, the UAV positions evolve continuously according to \eqref{eq:dynamics}, while the current belief is held fixed.
At the next measurement time, the candidate team configuration determines the EIG reward defined in \eqref{eq:cycle_eig_def}, and the incoming measurements update the belief used in the next cycle.

We account for continuation beyond that measurement through a repeated-cycle Bellman formulation.
Between measurements, the value function satisfies an HJB equation associated with the continuous UAV dynamics and a running cost.
At measurement time, a jump condition combines the EIG reward, a predicted covariance update, and the value of continuing from the post-measurement state.
A fixed point of this one-cycle Bellman map therefore represents repeated replanning without prescribing a finite terminal horizon.

To reduce the state-space dimensionality, we use a reduced approximation of a full belief-state control problem.
The current MAP estimates are treated as cycle-dependent parameters, while selected covariance descriptors are included in the numerical HJB state.
Each retained scalar covariance descriptor adds one dimension to the Cartesian HJB grid, which motivates the belief factorizations and parameter reductions introduced in Section~\ref{sec:eig_computation}.

The distinction between prediction and inference is essential: the fixed point propagates covariance summaries under the current nominal parameter estimates, whereas each actual measurement triggers a new joint estimate and a new HJB solve.
The following subsections define the running cost, the measurement update, and the approximate information reward before assembling the Bellman problem and the multi-agent feedbacks.

\subsection{Running Cost and Planning Objective}
\label{sec:running_cost}

Because UAV endurance is limited, energy expenditure is a practical bottleneck in prolonged search missions.
We therefore seek policies that localize the source efficiently not only in accuracy and time, but also in cumulative control effort relative to the ambient wind.
For a mission horizon $T>0$, we use the quadratic functional
\begin{equation}
\label{eq:control_effort}
E :=
\sum_{i=1}^M \int_0^T \norm{\u_i(t)-\bs w}^2\,\dd t.
\end{equation}
We interpret $E$ as a wind-relative control-effort proxy related to propulsion energy expenditure, not as a physical power or battery model.
It has units $\mathrm{m^2\,s^{-1}}$.
This proxy directly motivates the first term in the running cost below.
Over each interval between consecutive measurement times, the team accumulates the running cost
\begin{equation}
\label{eq:running_cost}
\begin{aligned}
L(\bX,\u,t)=
{}&\underbrace{\sum_{i=1}^M \frac{c_1}{2}\norm{\u_i-\bs w}^2}_{\text{wind-relative control effort}}
+\underbrace{\sum_{i=1}^M \frac{c_2}{2}\norm{\u_i-\u_p(\x_i,t)}^2}_{\text{patrol deviation}}
\\
&+\underbrace{c_3\!\!\sum_{1\le i<j\le M}\!\! f\!\left(\norm{\x_i-\x_j}\right)}_{\text{collision}}
+\underbrace{c_4\sum_{i,\ell} g_\ell(\x_i)}_{\text{no-fly zones}},
\end{aligned}
\end{equation}
% REVIEW C2-weights BEGIN
where $c_1,c_2>0$ and $c_3,c_4\ge0$ are user-specified weights; zero penalty weights disable the corresponding terms.
% REVIEW C2-weights END
We scale the objective to be dimensionless, so $L$ has units $\mathrm{s^{-1}}$: $c_1$ and $c_2$ have units $\mathrm{s\,m^{-2}}$, while $c_3$ and $c_4$ have units $\mathrm{s^{-1}}$ when the logistic potentials are dimensionless.
The coefficient $c_1$ scales the wind-relative control-effort proxy \eqref{eq:control_effort} relative to the other objective terms; it should not be interpreted as a calibrated conversion from velocity to physical energy.

The \emph{patrol field} $\u_p$ encourages systematic crosswind patrolling.
We use a sinusoidal path rotated to the wind direction,
\begin{equation}
\label{eq:patrol_field}
\u_p(\x,t)=
R(\alpha)\begin{bmatrix}-\sqrt{1-\rho^2}\\ \rho\sin(\omega t)\end{bmatrix}v_{\max},
\qquad
R(\alpha)=\begin{bmatrix}\cos\alpha&-\sin\alpha\\\sin\alpha&\cos\alpha\end{bmatrix},
\end{equation}
with $\rho\in(0,1)$ and $\omega=2\pi/T_p$.
This choice keeps $\norm{\u_p}\le v_{\max}$, gives a constant upwind component $\sqrt{1-\rho^{2}}v_{\max}$ opposite to $\bs w$, and superimposes a sinusoidal crosswind component of amplitude $\rho v_{\max}$.

The \emph{collision potential} $f$ is logistic with amplitude $\gamma_{\max}>0$:
\begin{equation}
\label{eq:collision_potential}
f(d)=\frac{\gamma_{\max}}{1+e^{ad-b}},\quad
f'(d)=-\gamma_{\max}a\frac{e^{ad-b}}{(1+e^{ad-b})^2}<0,
\end{equation}
% REVIEW C2-logistic BEGIN
calibrated from radii $0\le d_{\min}<d_{\max}$ and confidence $p_c\in(1/2,1)$:
% REVIEW C2-logistic END
\begin{equation}
\label{eq:collision_calibration}
a=\frac{2\ln(p_c/(1-p_c))}{d_{\max}-d_{\min}},\qquad
b=a\frac{d_{\max}+d_{\min}}{2}.
\end{equation}
The radii bracket the transition band: $f(d_{\min})=p_c\gamma_{\max}$ and $f(d_{\max})=(1-p_c)\gamma_{\max}$.
For $p_c$ close to one, the penalty is nearly saturated inside $d_{\min}$ and small beyond $d_{\max}$.
The no-fly penalties $g_\ell$ use the same logistic form, evaluated at a signed distance from the boundary of zone $\mathcal Z_\ell$.
Define $s_\ell(\x)=\operatorname{dist}(\x,\mathcal Z_\ell)$ for $\x\notin\mathcal Z_\ell$ and $s_\ell(\x)=-\operatorname{dist}(\x,\partial\mathcal Z_\ell)$ for $\x\in\mathcal Z_\ell$; then $g_\ell(\x)=\gamma_\ell/(1+\exp(a_\ell s_\ell(\x)-b_\ell))$, calibrated analogously to \eqref{eq:collision_calibration}.
The collision and no-fly terms are finite soft penalties and do not provide hard safety guarantees.

\subsection{Laplace Belief Update}
\label{sec:belief_update}

The full parameter vector is $\bs\xi=(\bs\theta,\zeta)\in\R^3$.
We use independent Gaussian prior factors $\mathcal N(\bs\mu_0,\bs\Sigma_0)$ for $\bs\theta$ and $\mathcal N(\bar\zeta_0,\bar\sigma_{\zeta,0}^2)$ for $\zeta$, restricted to the search box below.
Here the subscript $0$ on the prior parameters denotes their fixed values before the first observation; time-indexed quantities such as $\bs\Sigma_{t_0}$ denote the updated belief.

Because $C(\cdot,\bs\theta,\zeta)$ is nonlinear in $\bs\xi$, the posterior is non-Gaussian.
For the Gauss--Newton Laplace-type update~\cite{laplace_approx}, we find $\bs\xi_{\text{MAP}}=(\bs\theta_{\text{MAP}},\zeta_{\text{MAP}})$ by minimizing the negative log-posterior
\begin{equation}
\label{eq:neg_log_post}
\mathcal J(\bs\xi)=
\tfrac12(\bs\theta-\bs\mu_0)^\top\bs\Sigma_0^{-1}(\bs\theta-\bs\mu_0)
+\frac{(\zeta-\bar\zeta_0)^2}{2\bar\sigma_{\zeta,0}^2}
+\frac{1}{2\sigma_e^2}\sum_{i,j}
\!\left(y_{ij}-C(\x_i(t_j),\bs\theta,\zeta)\right)^2,
\end{equation}
where the sum includes all drones and measurement times in $D_k$.
We solve
\begin{equation}
\bs\xi_{\mathrm{MAP}}
\in \operatorname*{argmin}_{\bs\theta\in\Omega,\;\zeta\in[\zeta_{\min},\zeta_{\max}]}
\mathcal J(\bs\theta,\zeta).
\end{equation}
The prior normalizing constants are fixed and can be omitted from the constrained minimization.
At the MAP estimate, let $\bs J_C\in\R^{M(k+1)\times3}$ be the Jacobian whose rows are $\nabla_{\bs\xi}C(\x_i(t_j),\bs\xi_{\mathrm{MAP}})^\top$, and define the prior covariance $\bs P_0:=\diag(\bs\Sigma_0,\bar\sigma_{\zeta,0}^2)$.
The Gauss--Newton precision matrix is
\begin{equation}
\label{eq:gauss_newton_precision}
\bs H_{\mathrm{GN}}
=\bs P_0^{-1}+\frac{1}{\sigma_e^2}\bs J_C^\top\bs J_C,
\end{equation}
which retains the $\bs\theta$--$\zeta$ cross-terms from the likelihood.
Because $\bs P_0\succ0$ and $\bs J_C^\top\bs J_C\succeq0$, the matrix $\bs H_{\mathrm{GN}}$ is positive definite and therefore invertible.
The marginal covariance blocks are consequently
\begin{equation}
\label{eq:laplace_cov}
\bs\Sigma_{t_k} = \bigl(\bs H_{\mathrm{GN}}^{-1}\bigr)_{\bs\theta\bs\theta},
\qquad
\sigma_{\zeta,t_k}^2 = \bigl(\bs H_{\mathrm{GN}}^{-1}\bigr)_{\zeta\zeta}.
\end{equation}
The joint Gaussian covariance surrogate $\bs H_{\mathrm{GN}}^{-1}$ retains the location--intensity dependence.
It is evaluated at the constrained MAP estimate.
The Gaussian remains untruncated even when a bound is active, so its covariance omits the effects of the bounds on the posterior shape.
The marginal covariance blocks are passed to the planner.
The planner then uses a block-diagonal Fisher approximation to predict the effect of prospective future measurements (see \Cref{sec:eig_computation}).

\subsection{Linearized-Gaussian EIG and Planner Belief Jump}
\label{sec:eig_computation}

To evaluate the information objective \eqref{eq:cycle_eig_def} inside the control problem, the planner uses the factored Gaussian prior
\begin{equation}
\pi_t^-=\mathcal N(\bs\mu_t,\bs\Sigma_t)\times\mathcal N(\hat\zeta_t,\sigma_{\zeta,t}^2),
\end{equation}
and a first-order linearization of the measurement model.
This approximation yields a closed-form EIG and a predicted covariance update that can be evaluated without integrating over prospective measurement outcomes.

At the current nominal parameter
$\bs\xi^\star:=(\bs\mu_t,\hat\zeta_t)\in\R^3$, the exact
$3\times3$ Fisher information of the joint measurement model has cross-terms between $\bs\theta$ and $\zeta$,
\begin{equation}
\bs F^{3\times3}(\bX)
=\frac{1}{\sigma_e^2}\sum_{i=1}^M
\nabla_{\bs\xi}C(\x_i,\bs\xi^\star)\,
\nabla_{\bs\xi}C(\x_i,\bs\xi^\star)^\top,
\end{equation}
where $\bs\xi^\star$ is the current joint MAP estimate (cf.\ \eqref{eq:neg_log_post}).
A symmetric $3\times3$ joint covariance has six independent entries.
To avoid carrying all six covariance coordinates in the numerical control state, we adopt the following block-diagonal approximation:

\begin{itemize}
\item \textbf{(B1) Independent prior:} consistent with the marginalized planner prior $\pi_t^-$, the joint prior covariance is approximated by $\bs\Xi:=\diag(\bs\Sigma_t,\sigma_{\zeta,t}^2)$, i.e.\ $\bs\theta$ and $\zeta$ are treated as a priori independent at each planning cycle.
\item \textbf{(B2) Block-diagonal Fisher information:} the cross-terms $F_{\theta_j,\zeta}=\sigma_e^{-2}\sum_i(\partial_{\theta_j}C)(\partial_\zeta C)$ are dropped from $\bs F^{3\times3}$.
\end{itemize}

The resulting factored representation retains the three independent entries of the $2\times2$ location covariance and the intensity variance, but discards the two location--intensity cross-covariances.
It therefore reduces the number of scalar covariance coordinates from six to four.
After (B2) removes the cross-information terms, the surrogate precision is block-diagonal and the cycle-wise EIG factorizes:

\begin{proposition}[block-diagonal linearized-Gaussian planner EIG]
\label{prop:eig}
Under a first-order Taylor linearization of $C(\cdot,\bs\xi)$ about the current nominal parameter $\bs\xi^\star=(\bs\mu_t,\hat\zeta_t)$, additive Gaussian measurement noise, and the block-diagonal approximation (B1)--(B2), the cycle-wise EIG of the resulting Gaussian planner surrogate is
\begin{equation}
\label{eq:eig_formula}
\operatorname{EIG}(\bX;\pi_t^-)
=\underbrace{\tfrac12\log\det\!\bigl(\bs I_2+\bs\Sigma_t\,\bs F_{\bs\theta}(\bX)\bigr)}_{\text{location IG}}
+\underbrace{\tfrac12\log\bigl(1+\sigma_{\zeta,t}^2\,F_\zeta(\bX)\bigr)}_{\text{intensity IG}},
\end{equation}
where
\begin{align}
\bs F_{\bs\theta}(\bX)
&=\frac{1}{\sigma_e^2}\sum_{i=1}^M
\nabla_{\bs\theta}C(\x_i,\bs\mu_t,\hat\zeta_t)\,
\nabla_{\bs\theta}C(\x_i,\bs\mu_t,\hat\zeta_t)^\top
\;\in\R^{2\times2},
\label{eq:fisher_theta}\\
F_\zeta(\bX)
&=\frac{1}{\sigma_e^2}\sum_{i=1}^M
\Bigl(\tfrac{\partial C}{\partial\zeta}(\x_i,\bs\mu_t,\hat\zeta_t)\Bigr)^2
\;\in\R.
\label{eq:fisher_zeta}
\end{align}
\end{proposition}
\noindent A derivation for the location term is given in \Cref{app:eig}; the intensity term follows identically by treating $\zeta$ as a one-dimensional parameter.
For the plume model in \Cref{sec:plume_model}, $C(\x,\bs\theta,\zeta)=e^\zeta C(\x,\bs\theta,0)$, hence $\partial_\zeta C(\x,\bs\theta,\zeta)=C(\x,\bs\theta,\zeta)$.
Therefore $F_\zeta = \sigma_e^{-2}\sum_i C(\x_i,\bs\mu_t,\hat\zeta_t)^2$: the intensity Fisher information is the summed squared predicted signal divided by the noise variance, and the intensity EIG grows logarithmically with this quantity.

Retaining the $\bs\theta$--$\zeta$ cross-information terms, while still using the independent prior $\bs\Xi$ of (B1), gives the linearized-Gaussian expression
\begin{equation}
\operatorname{EIG}_{\mathrm{full}}
=\tfrac12\log\det\!\left(\bs I_3+\bs\Xi^{1/2}\bs F^{3\times3}\bs\Xi^{1/2}\right).
\end{equation}
\Cref{eq:eig_formula} is the planner surrogate used here after applying assumptions (B1)--(B2).
It has a definite directional bias.
Indeed, writing $\bs\Xi^{1/2}\bs F^{3\times3}\bs\Xi^{1/2}=\left[\begin{smallmatrix}\bs A&\bs C\\\bs C^\top&\bs B\end{smallmatrix}\right]\succeq0$, the Schur-complement identity gives
\begin{equation}
\det(\bs I_3+\bs\Xi^{1/2}\bs F^{3\times3}\bs\Xi^{1/2})
=\det(\bs I+\bs A)\det\!\left(\bs I+\bs B-\bs C^\top(\bs I+\bs A)^{-1}\bs C\right)
\le \det(\bs I+\bs A)\det(\bs I+\bs B).
\end{equation}
Thus, deleting the cross-information block makes the factored planner EIG an upper bound on $\operatorname{EIG}_{\mathrm{full}}$ and can overvalue configurations in which location and intensity are strongly confounded.
The comparison isolates (B2), since both sides use the same independent prior $\bs\Xi$; it does not bound the effect of (B1), because the Laplace belief $\bs H_{\mathrm{GN}}^{-1}$ generally has nonzero $\bs\theta$--$\zeta$ cross-covariances and no such ordering holds once the prior itself is correlated.

The factored form \eqref{eq:eig_formula} has two important consequences.
First, the closed form depends only on the candidate drone positions and current covariance summaries, so it can be evaluated without sampling future measurement realizations and inserted directly as the jump reward in the HJB.
Second, the belief state splits accordingly: under (B1)--(B2) the location and intensity uncertainties evolve independently at each jump,
\begin{equation}
\label{eq:belief_jump}
\bs\Sigma_+ = \bigl(\bs\Sigma_t^{-1}+\bs F_{\bs\theta}\bigr)^{-1},
\qquad
\sigma_{\zeta,+}^2 = \bigl(\sigma_{\zeta,t}^{-2}+F_\zeta\bigr)^{-1},
\end{equation}
so the planner can represent the intensity variance by one scalar covariance coordinate without coupling it to $\bs\Sigma$.
\Cref{eq:belief_jump} is the planner's local Gaussian/Fisher jump model used to predict the effect of a prospective measurement.
After a real measurement arrives, the actual belief is recomputed by the Laplace update \eqref{eq:neg_log_post}--\eqref{eq:laplace_cov} using the accumulated dataset.

\subsection{Repeated-Cycle HJB and Bellman Fixed Point}
\label{sec:hjb}

Since measurements arrive at the fixed interval $\Delta$, and the planner is re-invoked after every real observation, we consider a discounted one-cycle control problem conditioned on the current nominal belief means $(\bs\mu_t,\hat\zeta_t)$.
This should be interpreted as a tractable reduced approximation of the full belief-state control problem, rather than as propagation of the full conditional law between observation times.

\subsubsection{Control Formulation and Feedback}

% REVIEW A2-phase BEGIN
For a cycle beginning at $t_k$, we use the local time $\tau=t-t_k\in[0,\Delta]$ and write $L_k(\bX,\u,\tau):=L(\bX,\u,t_k+\tau)$.
The planning problem repeats this interval's patrol waveform in every virtual cycle, while actual replanning advances to the next physical-time interval.
The value function below belongs to the current problem $k$, whose index is suppressed; this waveform replay is an additional planning approximation.
% REVIEW A2-phase END
Let $r>0$ be an exponential discount rate, with units $\mathrm{s^{-1}}$, and let
\begin{equation}
S_\tau := (\bX_\tau,\bs\Sigma,\sigma_\zeta^2), \qquad \tau\in[0,\Delta),
\end{equation}
denote the reduced cycle-state process, where only the drone-position component $\bX_\tau$ evolves continuously between measurements while $\bs\Sigma$ and $\sigma_\zeta^2$ are held fixed over the cycle.
For a given cycle-start state $s=(\bX,\bs\Sigma,\sigma_\zeta^2)$, the post-measurement state is
\begin{equation}
S_\Delta^+ := (\bX_\Delta,\bs\Sigma_+,\sigma_{\zeta,+}^2),
\end{equation}
with $(\bs\Sigma_+,\sigma_{\zeta,+}^2)$ given by \eqref{eq:belief_jump}.
Let $V(\tau,s)$ denote the value at local cycle time $\tau$ and reduced state $s$.
Its cycle-start trace satisfies the self-referential Bellman equation
\begin{equation}
\label{eq:bellman}
V(0,s)
=
\inf_{\u\in\mathcal A}\E\!\left[
\int_0^\Delta e^{-r\tau}L_k(\bX_\tau,\u_\tau,\tau)\,\dd \tau
\;+\; e^{-r\Delta}\Bigl(V(0,S_\Delta^+) - \operatorname{EIG}(\bX_\Delta;\pi_{t_k}^-)\Bigr)
\;\middle|\; S_0=s
\right],
\end{equation}
where the expectation is taken over the drone-motion noise on the interval $[0,\Delta]$.
The EIG reward appears with a minus sign because the controller minimizes cost: informative measurements lower the objective, while motion cost and continuation also influence the chosen control.
The occurrence of $V(0,S_\Delta^+)$ on the right-hand side accounts for continuation after the next measurement.

On $0\le\tau<\Delta$, the same value function satisfies the Hamilton--Jacobi--Bellman PDE, interpreted in the viscosity sense~\cite{crandall1992viscosity},
\begin{equation}
\label{eq:hjb_pde}
-\partial_\tau V = \min_{\norm{\u_i}\le v_{\max}}
\Bigl[
L_k(\bX,\u,\tau)+(\nabla_{\bX}V)^\top\u
+\tfrac12\sigma_d^2\Delta_{\bX}V-rV
\Bigr],
\end{equation}
where
\begin{equation}
\Delta_{\bX}V:=\sum_{i=1}^M \Delta_{\x_i}V
=\sum_{i=1}^M\sum_{m=1}^2 \partial_{x_{i,m}x_{i,m}}^2 V
\end{equation}
is the Laplacian over the full team-position variable, and the diffusion term accounts for drone motion noise \eqref{eq:dynamics}.
The backward PDE is closed by the terminal jump condition
\begin{equation}
\label{eq:hjb_jump}
V(\Delta^-,s)
=V\bigl(0,\mathcal I(s)\bigr)
-\tfrac12\log\det\!\bigl(\bs I_2+\bs\Sigma\,\bs F_{\bs\theta}(\bX)\bigr)
-\tfrac12\log\bigl(1+\sigma_\zeta^2 F_\zeta(\bX)\bigr),
\end{equation}
where $\mathcal I(s)=(\bX,\bs\Sigma_+,\sigma_{\zeta,+}^2)$ is the deterministic planner jump defined by \eqref{eq:belief_jump}.
The factor $e^{-r\Delta}$ in \eqref{eq:bellman} is generated by the discount term $-rV$ during backward integration and therefore does not appear again in \eqref{eq:hjb_jump}.

Because $L$ is quadratic in each $\u_i$, the minimizer is available in closed form,
\begin{equation}
\label{eq:optimal_control}
\u_i^\star
=\Pi_{v_{\max}}\!\left(
\frac{c_1\bs w+c_2\u_p(\x_i,t_k+\tau)}{c_1+c_2}-\frac{\nabla_{\x_i}V(\tau,s)}{c_1+c_2}
\right),
\end{equation}
where $\Pi_{v_{\max}}$ denotes Euclidean projection onto the closed velocity ball,
\begin{equation}
\Pi_{v_{\max}}(a)=
\begin{cases}
a, & \|a\|\le v_{\max},\\[2mm]
\dfrac{v_{\max}}{\|a\|}a, & \|a\|>v_{\max}.
\end{cases}
\end{equation}

\subsubsection{Numerical Solution}

We pose the continuous spatial problem on $\R^{2M}$, with bounded-value behavior inherited from the bounded controls, running cost, and jump reward.
% REVIEW B5-box BEGIN
Numerically, the relative spatial box is $D_z=(\Omega-\bs\mu_t)^M$: each coordinate endpoint is the corresponding physical-domain endpoint minus the current mean.
Thus a drone inside $\Omega$ lies inside its relative spatial grid.
Absolute no-fly penalties are evaluated at $\x_i=\z_i+\bs\mu_t$.
Off-grid interpolation foot points are clipped to the computational box; this numerical closure is distinct from projecting simulated physical positions into $\Omega$ and from clamping covariance coordinates.
% REVIEW B5-box END

For a fixed cycle $k$, let $\mathcal B_k$ denote the map that takes a candidate cycle-start trace to the new time-zero trace produced by one backward solution of \eqref{eq:hjb_pde} with the jump \eqref{eq:hjb_jump}.
We again suppress $k$ when discussing one such solve.
% REVIEW B1-contraction BEGIN
For the continuous reduced control problem with well-defined transitions and bounded running cost and jump reward, the Bellman map is a sup-norm contraction with factor $e^{-r\Delta}<1$, so its Picard iteration
% REVIEW B1-contraction END
\begin{equation}
\label{eq:picard}
V_{n+1} = \mathcal B[V_n],\qquad V_0=0,\qquad n=0,1,2,\ldots
\end{equation}
converges geometrically to the cycle-start trace $V(0,\cdot)$.
Here $V_n$ approximates the cycle-start trace, and $V_0=0$ initializes the fixed-point iteration.
The backward PDE uses the jump condition \eqref{eq:hjb_jump} at the end of each cycle.
% REVIEW B1-discrete BEGIN
Each Bellman application requires one backward HJB solve on a Cartesian grid over $[0,\Delta]$.
For a numerical step of length $h$, let $W$ be the next-time array, $\mathcal Q_h$ clipped multilinear interpolation, $D_h$ the spatial discrete Laplacian, and $\u_h$ the control \eqref{eq:optimal_control} evaluated using finite-difference gradients of $W$.
The implemented update is
\begin{equation}
\label{eq:discrete_hjb_update}
W^- = e^{-rh}\mathcal Q_h[W](s+h(\u_h,0))
+\frac{1-e^{-rh}}{r}\left(L_k(\bX,\u_h,\tau)+\tfrac12\sigma_d^2D_hW\right).
\end{equation}
Belief coordinates remain fixed in this step, and the control is selected from the value gradient.
The continuous contraction result does not directly establish contraction of this discrete update; its residual is used as an iteration stopping criterion.
The finite-difference and feedback-evaluation conventions are specified in \Cref{sec:computational_settings}.
% REVIEW B1-discrete END
In practice, Picard iteration can be replaced by Anderson acceleration applied to the residuals $R_n:=\mathcal B[V_n]-V_n$, using a short window of previous iterates to extrapolate the next update while targeting the same Bellman fixed point~\cite{anderson1965iterative}.
That is, after computing several residuals $R_n$, the next iterate is formed from a small least-squares combination of recent Bellman images rather than from the raw Picard step alone.
Anderson acceleration with depth $m=3$ is used in all reported experiments (see \Cref{sec:setup}).
The implemented acceptance statistic is the normalized Bellman residual
\begin{equation}
\label{eq:normalized_bellman_residual}
\rho_n:=\frac{\norm{\mathcal B[V_n]-V_n}_\infty}{1+\norm{V_n}_\infty}.
\end{equation}
A solve is accepted when $\rho_n<10^{-3}$; otherwise it is retried from scratch with $\Delta t_{\mathrm{HJB}}/2$ and the same iteration cap.
% REVIEW C1-retry BEGIN
An Anderson proposal is accepted only if it is finite and does not increase $\rho_n$ relative to the raw Picard proposal.
If the retry also fails, the implementation records nonconvergence and returns its final iterate; such cases are not accepted in the reported comparisons.

% REVIEW C1-retry END

After each real measurement at $t_k$, the Laplace update \eqref{eq:neg_log_post} produces updated belief parameters $(\bs\mu_{t_k},\bs\Sigma_{t_k},\sigma_{\zeta,t_k}^2)$.
The fixed-point \eqref{eq:picard} is then resolved from scratch with the new belief, and the resulting policy is applied until $t_{k+1}$.

\subsubsection{State Representation}

During one HJB solve, covariance descriptors vary across the grid and change through the predicted measurement jumps; they remain constant along each between-measurement trajectory.
The observed covariances determine the cycle-start state and grid extent, while the current parameter estimates, wind, plume parameters, and cost weights remain fixed throughout the solve.
Positions are represented relative to the current location estimate, $\z_i=\x_i-\bs\mu_t$.
The covariance $\bs\Sigma$ is encoded by its log-eigenvalues $(\log\lambda_1,\log\lambda_2)$ and principal angle $\phi$; together with the log-intensity variance $\log\sigma_\zeta^2$, this yields the reduced numerical HJB state
\begin{equation}
s=\bigl(\z_1,\ldots,\z_M,\,\log\lambda_1,\log\lambda_2,\phi,\,\log\sigma_\zeta^2\bigr)
\in\R^{2M+4}.
\end{equation}
% REVIEW C1-coordinates BEGIN
At the measured belief, $\lambda_1\ge\lambda_2>0$ are its ordered eigenvalues and $\phi$ is the principal-eigenvector angle relative to the $x_1$-axis, understood modulo $\pi$.
These descriptors define the covariance coordinates; their numerical grid and interpolation conventions are given in \Cref{sec:computational_settings}.
% REVIEW C1-coordinates END
We refer to the four scalar covariance descriptors
$\log\lambda_1$, $\log\lambda_2$, $\phi$, and $\log\sigma_\zeta^2$ as the \emph{belief axes} of the Cartesian HJB grid.

\subsection{Scalability Approximations}
\label{sec:scalability}

The $2M+4$-dimensional HJB state of the joint $M$-agent factored planner is tractable only for small $M$.
Rather than solving this joint problem, we introduce two complementary approximations that restrict the solve to low-dimensional subproblems independent of $M$.

\subsubsection{Wind-Aligned Belief Parameterization}
\label{sec:wind_aligned}
\begin{figure}[htbp]
  \centering
  \includegraphics[width=.8\linewidth]{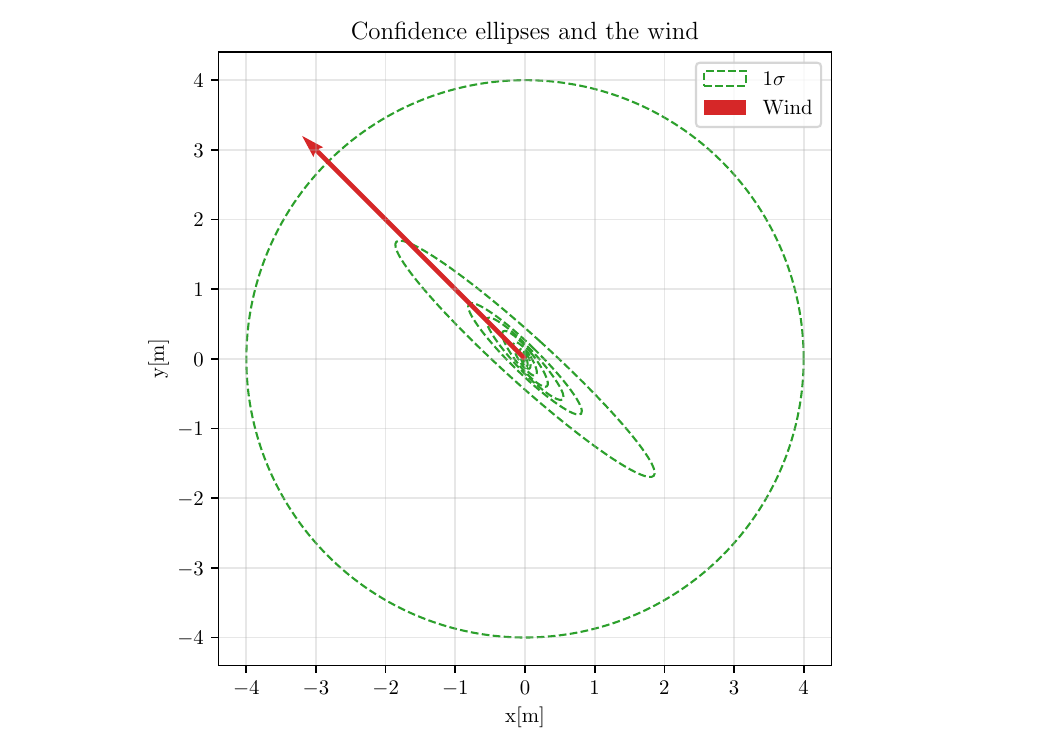}
  \caption{Successive $1\sigma$ belief contours after repeated observations, together with an arrow indicating the wind direction.
  The dominant direction of uncertainty is visibly aligned with the wind.
  This is exploited to reduce the dimensionality of the HJB.}
  \label{fig:wind_as_angle}
\end{figure}

We parameterize the posterior covariance by its two eigenvalues and its principal-axis angle $\phi$.
A natural approximation, illustrated in \Cref{fig:wind_as_angle}, is to fix $\phi=\alpha$ (the plume-axis direction), removing the angle axis and reducing state dimension from $2M+4$ to $2M+3$.
This eliminates one grid dimension from the HJB, saving a factor of $n_b$ in compute (the number of grid points along that axis).
% REVIEW A3-mapping BEGIN
Initialization uses the measured covariance's ordered eigenvalues with the fixed wind orientation.
At each predicted jump, the planner decodes the covariance, forms the Fisher blocks, applies \eqref{eq:belief_jump}, extracts the ordered updated eigenvalues, and reorients them onto the wind axes.
This spectral reorientation differs from deleting off-diagonal entries in the wind basis.
The resulting covariance is then encoded for continuation interpolation using the grid conventions in \Cref{sec:computational_settings}.
The alignment is an empirical approximation; \Cref{sec:exp_wind_alignment} measures its effect on localization and planning time.
% REVIEW A3-mapping END

\subsubsection{One-Body and Pairwise Planning Surrogates}
\label{sec:scalability_surrogates}
Rather than solving the full $M$-agent HJB, we define two low-dimensional surrogate policies whose per-cycle HJB subproblem dimension is independent of team size.
Both policies control the full team using a shared belief.
One-body planning evaluates each drone's prospective measurements independently within the HJB subproblem.
Pairwise planning evaluates two drones together, allowing their combined information and interaction costs to influence the feedback.
Aggregating the pairwise value gradients then gives each drone a command that depends on the positions of the other team members.

\emph{One-body planning.}
A single $(2+d_b)$-dimensional fixed-point problem is solved, with $d_b=3$ accounting for the two log-eigenvalues and the log-intensity variance (the principal angle is fixed by the wind-aligned approximation).
From its solution we obtain a feedback $\u^{(1)}(\z_i,\bs\Sigma,\sigma_\zeta^2)$ that is reused for every drone independently.
Because the one-body HJB describes only a single drone, pairwise collision structure is absent; it must be reintroduced via an explicit gradient correction,
\begin{equation}
\label{eq:one_body_collision_correction}
% REVIEW A1-gain: added collision time scale
\u_i \leftarrow \Pi_{v_{\max}}\!\left(
\u_i^{(1)}(\z_i,\bs\Sigma,\sigma_\zeta^2)
-\tau_{\mathrm{col}}\frac{c_3}{c_1+c_2}\nabla_{\x_i}\!\sum_{j\ne i}f(\norm{\x_i-\x_j})
\right),
\end{equation}
% REVIEW A1-convention BEGIN
where $\tau_{\mathrm{col}}>0$ is a time scale, independent of the simulation and HJB timesteps, making the correction a velocity.
The nominal feedback is projected first, then corrected and projected again; the minus sign gives repulsion because $f'(d)<0$.
For numerical evaluation, the pair gradient uses $f'(r_\epsilon)(\x_i-\x_j)/r_\epsilon$, with $r_\epsilon=\max(\norm{\x_i-\x_j},10^{-8}\,\mathrm{m})$.
The correction vanishes at exact coincidence, and the safeguard does not guarantee separation.
% REVIEW A1-convention END

\emph{Pairwise planning.}
A $(4+d_b)$-dimensional two-agent HJB problem is solved, yielding a pair feedback $\u^{(2)}(\z_i,\z_j,\bs\Sigma,\sigma_\zeta^2)\in\R^{4}$.
Let $V^{(2)}$ denote the corresponding pairwise value function and define the aggregated value gradient for drone $i$ by
\begin{equation}
\label{eq:pairwise_control}
\bar g_i
:=\frac{1}{M-1}\sum_{j\ne i}
\nabla_{\z_i}V^{(2)}(\z_i,\z_j,\bs\Sigma,\sigma_\zeta^2),
\qquad
\u_i
=\Pi_{v_{\max}}\!\left(
\frac{c_1\bs w+c_2\u_p(\x_i,t)}{c_1+c_2}
-\frac{\bar g_i}{c_1+c_2}
\right).
\end{equation}
Thus, the implementation averages pairwise value gradients and applies the quadratic Hamiltonian minimizer and velocity projection once; it does not average already projected pair controls.
Collision and no-fly penalties enter $\bar g_i$ through the two-agent HJB, and the final multi-agent command uses the same runtime safeguard \eqref{eq:one_body_collision_correction}, with $\tau_{\mathrm{col}}=1\,\mathrm{s}$ for both controllers in the reported experiments.

\Cref{eq:pairwise_control} is a heuristic aggregation of two-agent feedbacks, not the exact minimizer of a globally defined $M$-agent HJB.
The one-body and pairwise policies are fixed-dimensional feedback surrogates whose empirical behavior is evaluated in \Cref{sec:experiments}; they are not error-controlled approximations of the full $M$-agent value function.

For both surrogates, the per-cycle HJB subproblem dimension is independent of $M$: the one-body subproblem is always $(2+d_b)$-dimensional (5-dimensional in the experiments reported here) and the pairwise subproblem is always $(4+d_b)$-dimensional (7-dimensional), regardless of team size.
The HJB table dimension and the cost of constructing one table are independent of $M$.
At online evaluation, one-body table queries cost $O(M)$ and pairwise table queries cost $O(M^2)$.
The all-pairs collision correction makes evaluation of either complete command $O(M^2)$.
These costs are separate from joint data assimilation, whose dataset grows with team size and mission duration.

\subsection{Closed-Loop Procedure}
\label{sec:closed_loop_procedure}

\Cref{alg:closed_loop} summarizes the complete closed-loop repeated-cycle HJB controller.
\begin{algorithm}[t]
\caption{Closed-loop repeated-cycle HJB controller}
\label{alg:closed_loop}
\begin{algorithmic}[1]
\State \textbf{Initialize} belief parameters $(\bs\mu_0,\bs\Sigma_0,\hat\zeta_0,\sigma_{\zeta,0}^2)$, drone positions $\{\x_i(0)\}_{i=1}^M$, and accumulated dataset $D_{-1}=\emptyset$.
% REVIEW C1-B4-algorithm BEGIN
\If{the initial time precedes the first scheduled observation}
    \State Solve using the prior and propagate feedback until that observation or $T$.
\EndIf
\For{each scheduled measurement time $t_k\le T$, including the initial time if configured}
    \State Collect $y_{ik}$; compute \eqref{eq:stop_residual} with the pre-update estimates; form $D_k$.
    \State Solve \eqref{eq:neg_log_post} and form the updated covariance via \eqref{eq:laplace_cov}.
    \If{the full stopping rule in \Cref{sec:stopping_criterion} is met, or $t_k=T$}
        \State \textbf{break}
    \EndIf
    \State Build the relative grid; fix nominal means, wind, plume, costs, and this cycle's patrol waveform.
    \State Solve \eqref{eq:bellman}--\eqref{eq:picard} with the documented retry convention.
    \State Propagate feedback to the next observation or $T$.
\EndFor

% REVIEW C1-B4-algorithm END
\end{algorithmic}
\end{algorithm}

\section{Experiments}
\label{sec:experiments}

The experiments examine the computational and mission-level trade-offs introduced by the reduced controller.
We first measure the cost and accuracy of wind alignment, then compare one-body HJB control with patrol and greedy EIG policies.
The remaining studies test whether pairwise feedback improves cooperative localization, how discounting changes the accuracy--effort trade-off, and how the pairwise controller behaves as the team grows from two to four drones.

\subsection{Experimental Setup}
\label{sec:setup}

\subsubsection{Simulation Scenarios}
\label{sec:simulation_scenarios}
\begin{figure}[htbp]
  \centering
  \includegraphics[width=0.72\linewidth]{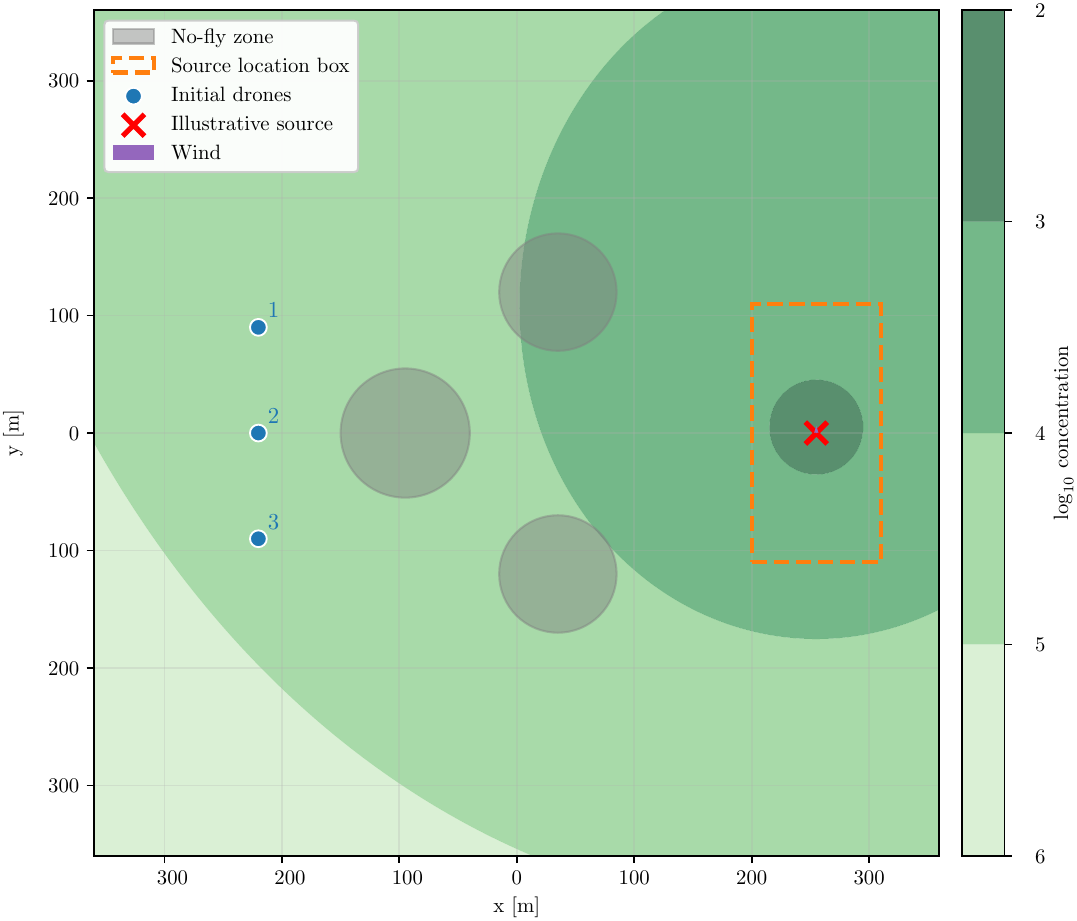}
  \caption{Challenging cooperative scenario geometry.
  The figure shows the initial drone positions, the source-spawn box, the no-fly zones, the weak vertical wind, and plume contours for an illustrative source at the center of the source box.}
  \label{fig:pairwise_showcase_setting}
\end{figure}

We use two simulation scenarios.
The first is a baseline scenario designed to compare one-body HJB replanning against non-adaptive patrol and greedy one-step EIG.
The second is a more challenging cooperative scenario designed to test pairwise planning and discount sensitivity and is illustrated in \Cref{fig:pairwise_showcase_setting}.
Both scenarios use three drones unless stated otherwise, the sliced screened plume model from \Cref{eq:plume_model}, joint source-location/log-intensity inference, Anderson-accelerated fixed-point iteration, and the diagnostic stopping structure described in \Cref{sec:stopping_criterion}.
% REVIEW C2-sampling BEGIN
Their numerical parameters are summarized in \Cref{tab:simulation_scenarios}.
In each scenario the source coordinates are sampled independently and uniformly over the stated source box.
% REVIEW C2-sampling END

\begin{table}[p]
  \centering
  \footnotesize
  \setlength{\tabcolsep}{3pt}
  \renewcommand{\arraystretch}{1.14}
  \caption{Simulation scenarios used in the numerical experiments.
  Units are included in the quantity column; shared entries are merged across the two scenario columns.}
  \label{tab:simulation_scenarios}
  \begin{tabular}{@{}>{\raggedright\arraybackslash}p{0.23\linewidth}>{\raggedright\arraybackslash}p{0.34\linewidth}>{\raggedright\arraybackslash}p{0.34\linewidth}@{}}
    \toprule
    Quantity & Baseline & Challenging cooperative \\
    \midrule
    Purpose & Baseline comparisons (one-body HJB, greedy EIG, patrol) & One-body vs.\ pairwise planning; discount sensitivity \\
    Drones & $M=3$ & $M=3$ \\
    Domain coordinates [m] % REVIEW C2-units
 & $[-500,500]^2$ & $[-360,360]^2$ \\
    Initial drone positions [m] & $\begin{array}{@{}l@{}}(-450,400),\,(-400,-100),\\(-450,-400)\end{array}$ & $\begin{array}{@{}l@{}}(-220,90),\,(-220,0),\\(-220,-90)\end{array}$ \\
    Source-box coordinates [m] % REVIEW C2-units
 & $[150,400]\times[-250,250]$ & $[200,310]\times[-110,110]$ \\
    Source intensity & chosen to match $\mathrm{SNR}_0\in\{50,100,150,200\}$ & $q=60$ \\
    Prior mean [m] & $(0,0)$ & $(255,0)$ \\
    Prior covariance [$\mathrm{m}^2$] & $300^2 I_2$ & $\diag(90^2,110^2)$ \\
    Log-intensity prior and bounds & $\mathcal N(\log 30,5.0)$ on $[2,10]$ & $\mathcal N(\log 45,1.0)$ on $[\log 10,\log 500]$ \\
    Wind [$\mathrm{m\,s^{-1}}$] & $(-5,1.5)$ & $(0,0.5)$ \\
    Plume parameters & $d_s=80\,\mathrm{m^2\,s^{-1}},\ \tau_s=3000\,\mathrm{s},\ \varepsilon_c=5\,\mathrm{m}$ & $d_s=90\,\mathrm{m^2\,s^{-1}},\ \tau_s=2600\,\mathrm{s},\ \varepsilon_c=22\,\mathrm{m}$ \\
    Measurement noise std. & $\sigma_e=5\times10^{-5}$ & $\sigma_e=2\times10^{-5}$ \\
    No-fly zones [m] & none & circles $(-95,0),R=55$, $(35,120),R=50$, $(35,-120),R=50$ \\
    Collision model [m] & $d_{\min}=50,\ d_{\max}=70,\ \gamma_{\max}=100,\ p_c=0.99$ & $d_{\min}=45,\ d_{\max}=90,\ \gamma_{\max}=80,\ p_c=0.99$ \\
    Cost weights & $c_1=c_2=10^{-4},\ c_3=10^{-2},\ c_4=0$ & $c_1=10^{-4},\ c_2=10^{-5},\ c_3=10^{-2},\ c_4=1.8$ \\
    Collision time scale & \multicolumn{2}{c}{$\tau_{\mathrm{col}}=1\,\mathrm{s}$} \\ % REVIEW A1-value
    Patrol parameters & $T_p=200\,\mathrm{s},\ \rho=0.8$ & $T_p=220\,\mathrm{s},\ \rho=0.82$ \\
    Velocity bound [$\mathrm{m\,s^{-1}}$] & $v_{\max}=12$ & $v_{\max}=10$ \\
    Measurement interval [s] & $\Delta=5$ & $\Delta=30$ \\
    Simulation timestep [s] & $\Delta t=0.05$ & $\Delta t=0.1$ \\
    Horizon [s] & $T=400$ & $T=600$ \\
    Dynamics noise std. [$\mathrm{m\,s^{-1/2}}$] & $\sigma_d=0.5$ & $\sigma_d=0.3$ \\
    Discount rate [$\mathrm{s^{-1}}$] & $r=0.3\,\mathrm{s^{-1}}$ & $r=0.3\,\mathrm{s^{-1}}$ for pairwise showcase; swept in discount study \\
    Fixed-point solver & Anderson $m=3$, tolerance $10^{-3}$, maximum $40$ iterations & Anderson $m=3$, tolerance $10^{-3}$, maximum $12$ iterations \\
    HJB discretization & nominal $\Delta x=50\,\mathrm{m}$ (20 nodes; actual $52.63\,\mathrm{m}$), nominal coarse spacing $100\,\mathrm{m}$ (10 nodes; actual $111.11\,\mathrm{m}$), $\Delta t_{\mathrm{HJB}}=0.1\,\mathrm{s}$ & nominal $\Delta x=60\,\mathrm{m}$ (12 nodes; actual $65.45\,\mathrm{m}$), nominal coarse spacing $120\,\mathrm{m}$ (6 nodes; actual $144\,\mathrm{m}$), $\Delta t_{\mathrm{HJB}}=1\,\mathrm{s}$, $3$ belief grid points per uncertainty axis \\
    \bottomrule
  \end{tabular}
\end{table}

In the baseline scenario, the source intensity $q$ is chosen analytically so that the \emph{source-point signal-to-noise ratio}
\begin{equation}
\label{eq:snr0_definition}
\mathrm{SNR}_0
:=
\frac{C(\bs\theta,\bs\theta,\zeta)}{\sigma_e}
\end{equation}
matches the requested level.
For the sliced screened plume \eqref{eq:plume_model}, evaluating at $\x=\bs\theta$ gives
\begin{equation}
\label{eq:snr0_sliced_plume}
C(\bs\theta,\bs\theta,\zeta)
=
\frac{e^\zeta}{4\pi d_s \varepsilon_c}\exp(-\kappa\varepsilon_c),
\end{equation}
so $\mathrm{SNR}_0$ is linear in $q=q_{\mathrm{ref}}e^\zeta$ and can be inverted directly when constructing each randomized scenario.

The challenging cooperative scenario places the source away from the initial drone positions and introduces obstacles along the search region.
It tests whether pairwise planning and weaker discounting help the team reach informative measurement locations while managing motion cost and operational penalties.

\subsubsection{Outcome Metrics and Stopping Criterion}
\label{sec:stopping_criterion}

To measure control effort, we approximate the continuous functional $E$ from \eqref{eq:control_effort} by the discrete quantity $\mathcal E$.
The drone dynamics are integrated over $N_t$ time steps of size $\Delta t$, and the velocity command $\u_i^n$ is held fixed over the $n$-th simulation step, giving the Riemann sum
\begin{equation}
\mathcal{E} = \sum_{i=1}^M \sum_{n=0}^{N_t-1} \norm{\u_i^n-\bs w}^2\,\Delta t.
\end{equation}
% REVIEW C2-distance BEGIN
Thus, $\mathcal E$ is the numerical approximation reported for $E$, and both quantities have units $\mathrm{m^2\,s^{-1}}$.
Reported travel distance is the sum $\sum_{i,n}\norm{\x_i^{n+1}-\x_i^n}$ over simulated position increments, including noise and physical-domain projection, at the fixed timestep in \Cref{tab:simulation_scenarios}.
% REVIEW C2-distance END

Early termination in the numerical experiments is not triggered solely by the location covariance becoming small.
The stopping rule combines the residual of the new batch, predicted using the pre-assimilation estimates at $t_k^-$, with the location and intensity variances after assimilation.
Its three conditions are
\begin{equation}
\label{eq:stop_sigma}
\tr(\bs\Sigma_{t_k}) < \epsilon_\sigma,
\end{equation}
\begin{equation}
\label{eq:stop_residual}
% REVIEW B4-residual: evaluate at pre-assimilation estimates
\frac{1}{M}\sum_{i=1}^M
\left(
\frac{y_{ik}-C(\x_i(t_k),\bs\mu_{t_k^-},\hat\zeta_{t_k^-})}{\sigma_e}
\right)^2
< \epsilon_{\mathrm{res}},
\end{equation}
\begin{equation}
\label{eq:stop_zeta}
\sigma_{\zeta,t_k}^2 < \epsilon_\zeta.
\end{equation}
Each condition serves a distinct purpose.
The trace criterion \eqref{eq:stop_sigma} measures whether the source-location belief has contracted enough for the localization task itself.
% REVIEW B4-timing BEGIN
The residual criterion \eqref{eq:stop_residual} checks the new batch against the predicted concentrations.
% REVIEW B4-timing END
The intensity-variance criterion \eqref{eq:stop_zeta} ensures that the joint inference problem has not left the source strength essentially undetermined even though the location covariance alone appears small.
Together, these three checks define a heuristic diagnostic stopping rule that is more restrictive than a location-covariance threshold alone.

The baseline scenario uses
\begin{equation}
\epsilon_\sigma=10\,\mathrm{m}^2,\qquad
\epsilon_{\mathrm{res}}=3,\qquad
\epsilon_\zeta=0.5,
\end{equation}
after at least $20$ measurement updates and for $3$ consecutive satisfied checks.
The challenging cooperative scenario uses the less restrictive thresholds
\begin{equation}
\epsilon_\sigma=40\,\mathrm{m}^2,\qquad
\epsilon_{\mathrm{res}}=8,\qquad
\epsilon_\zeta=2,
\end{equation}
after at least $6$ measurement updates and for $2$ consecutive satisfied checks.

In the experiments below, a run may terminate early if this diagnostic criterion is met; otherwise it continues until the fixed horizon $T$.
In the baseline comparison, a trial is successful if the final source-location estimate $\bs\mu_{\mathrm{final}}$ satisfies
\begin{equation}
\|\bs\mu_{\mathrm{final}}-\bs\theta\|\le 50\,\mathrm{m},
\end{equation}
and unsuccessful otherwise.
For the cooperative studies, we report final localization error directly together with success rates at several tighter error thresholds.

\subsubsection{Statistical Protocol}

All randomized comparisons use a paired seed protocol whenever multiple methods are evaluated on the same scenario: a common seed fixes the source location and measurement-noise realization, so differences between methods are due to the feedback policy.
For a given seed, source sampling and the simulation-noise stream are initialized identically across methods; the deterministic MAP-start generator is indexed by the same seed and measurement time.
Thus, paired methods share source draws, noise draws, and MAP candidate sets, although their different trajectories evaluate the common noise sequence at different positions.
For binomial rates, we report two-sided Wilson score $95\%$ confidence intervals \cite{wilson1927probable}.
For medians and medians of paired seed-level differences, we report two-sided bias-corrected and accelerated (BCa) bootstrap $95\%$ confidence intervals \cite{efron1987better}, using $10{,}000$ resamples of the independent seeds.
The shaded $10$th--$90$th percentile bands in the figures show distributional spread across trials and should not be interpreted as inferential confidence intervals.

\subsubsection{Computational Settings}
\label{sec:computational_settings}

The MAP problem is solved by constrained multi-start quasi-Newton optimization (L-BFGS-B) over the source-location and log-intensity bounds.
The optimizer is initialized from 256 uniformly sampled candidate source locations in $\Omega$, augmented by the current belief mean and the initial prior mean, with the log-intensity always initialized from the current estimate.
The negative log-posterior is evaluated at all candidates; only the 8 best are refined by the local optimizer, and the lowest resulting optimum is retained.
% REVIEW C1-MAP-failure BEGIN
If all top candidates yield non-finite results, the single best initial guess is refined as a fallback.
Optimizer status is recorded; a finite best result or the fallback's returned estimate is used even when the optimizer reports failure.
% REVIEW C1-MAP-failure END

At every real measurement time the Bellman fixed point is recomputed from scratch with the updated belief; no value table is reused across measurement cycles or parameter values.
% REVIEW B5-grid-extent BEGIN
The HJB grid is built over the translated spatial box $(\Omega-\bs\mu_t)^M$ described in \Cref{sec:hjb}.
% REVIEW B5-grid-extent END
For multi-drone HJB grids, the first drone's spatial axes use the primary node count, while all other drone spatial axes use the listed coarse node count.
The log-eigenvalue and log-intensity-variance axes each span a factor of $1.5$ below their current value, that is, $[\log(\lambda/1.5),\log\lambda]$ for a covariance eigenvalue $\lambda$ and likewise for $\log\sigma_\zeta^2$.

% REVIEW A3-grid-closure BEGIN
Predicted measurement jumps that leave a log-variance interval are clamped to its endpoints before multilinear interpolation of the continuation value.
The information reward still uses the untruncated Fisher update, so it need not equal the entropy reduction of the represented covariance transition; positive reward may persist at a variance floor.
This is a finite-grid approximation to continuation, separate from the actual Laplace update after each observation.
% REVIEW A3-grid-closure END

% REVIEW C1-coordinate-grid BEGIN
The two variance coordinates are decoded directly at Cartesian grid nodes, including combinations with reversed ordering, and sorted again after a predicted measurement update.
In the full-angle comparison, the angle grid has the same node count as a belief axis and spans $[\phi_0-\pi/4,\phi_0+\pi/4]$, where $\phi_0$ is the measured principal angle.
Updated angles are wrapped modulo $\pi$ from the lower endpoint and clamped to this local interval.
% REVIEW C1-coordinate-grid END

% REVIEW C1-evaluation BEGIN
Spatial derivatives use centered interior and one-sided boundary differences; $D_h$ applies the first-derivative operation twice along each spatial axis, using second-order edge formulas when at least three nodes are available.
Spatial interpolation foot points are clipped at the computational boundary.
One-body feedback differentiates a space--time multilinear value interpolant with perturbations of $0.01$ grid spacings and clipped endpoints; pairwise feedback interpolates stored spatial gradients in space and time before aggregation.
Residual evaluation and the final backward solve require additional Bellman applications beyond the reported solver iterations.
% REVIEW C1-evaluation END

% REVIEW B3-configs BEGIN
The archived run configurations, trajectories, and reproduction instructions are available in the v0.1.0 computational archive
\footnote{
\href{https://github.com/agcarlon/drone-soc-paper/tree/v0.1.0}{https://github.com/agcarlon/drone-soc-paper/tree/v0.1.0}.
}.
The reported results were generated with revision \texttt{f1ff89c}, whose source is included in the archive alongside the revised controller and regression tests for the explicit collision time scale.
The computations were performed on an Apple M2 Max workstation with $12$ physical and logical CPU cores and $64\,\mathrm{GiB}$ RAM, running macOS on arm64.
% REVIEW B3-configs END

\subsection{Wind-Aligned Belief Parameterization: Computational Savings}
\label{sec:exp_wind_alignment}

The wind-aligned approximation $\phi=\alpha$ (\Cref{sec:wind_aligned}) fixes the principal axis of the posterior covariance to the plume direction, eliminating $\phi$ as a state variable.
For one-body planning this reduces the HJB grid from 6-D to 5-D;
for pairwise planning, from 8-D to 7-D.
Since the fixed-point iteration scales linearly in the grid volume, every eliminated dimension saves a factor of $n_b$ operations (the number of grid points along the belief axes, typically $3$--$5$).

The approximation is motivated empirically, not analytically.
The Fisher information geometry depends on the drone--source geometry, so the reduction may affect both planning time and mission performance.
To verify this, we run one-body planning on $100$ seeds of the challenging cooperative scenario and compare two belief parameterizations: a 5-D wind-aligned HJB solve, in which $\phi$ is fixed to $\alpha$, and a 6-D full-angle HJB solve, in which $\phi$ is retained as an additional belief coordinate.
All other parameters are identical, including three grid points per belief uncertainty axis and $\Delta=30\,\mathrm{s}$.

\begin{figure}[htbp]
  \centering
  \includegraphics[width=\linewidth]{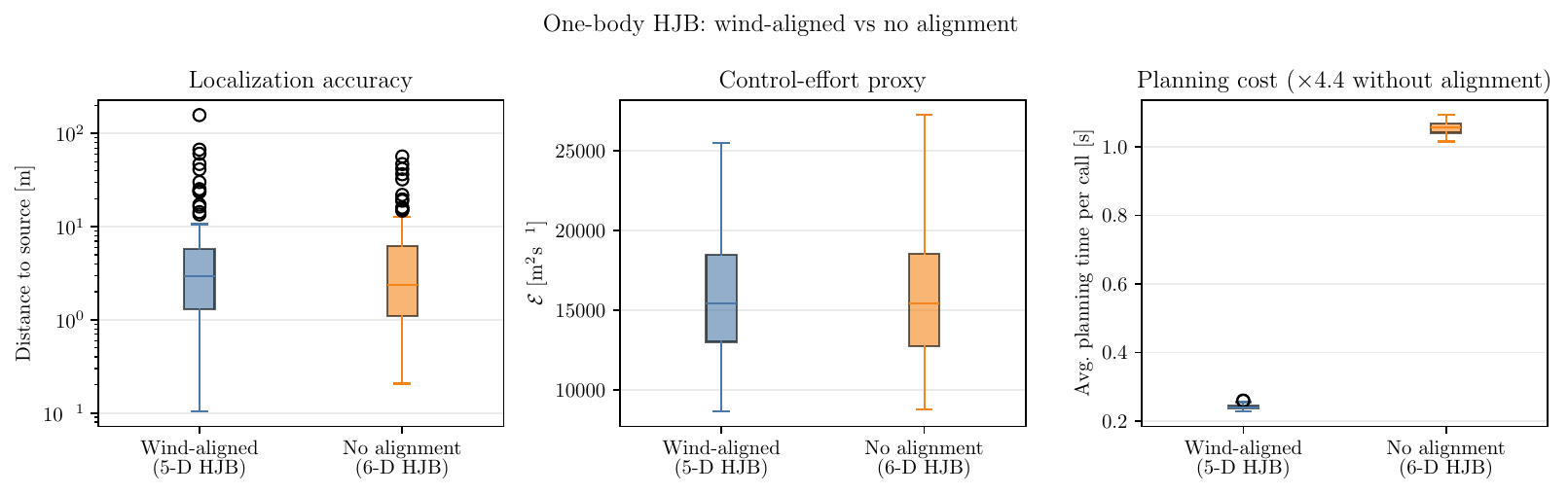}
  \caption{One-body planning: wind-aligned ($5$-D HJB, 3{,}888 grid points) versus no alignment ($6$-D HJB, 11{,}664 grid points), across $100$ seeds of the challenging cooperative scenario.
  Left and center: wind alignment gives a slightly larger median distance to the source, while the control-effort proxy $\mathcal{E}$ is nearly unchanged.
  Right: the median planning time per call is $4.4\times$ smaller with wind alignment.
  Horizontal bars indicate medians.}
  \label{fig:wind_aligned_comparison}
\end{figure}

Results across $100$ seeds are shown in \Cref{fig:wind_aligned_comparison}.
The median distance to the source is $2.92\,\mathrm{m}$ with wind alignment and $2.37\,\mathrm{m}$ without, while the diagnostic-stop rates are $85\%$ and $91\%$, respectively.
The median control efforts are nearly equal, $1.54\times10^4\,\mathrm{m^2s^{-1}}$ for both parameterizations.
By contrast, the median planning time per call is $0.24\,\mathrm{s}$ with wind alignment versus $1.06\,\mathrm{s}$ without, a $4.4\times$ speedup.
The grid-volume ratio is exactly $n_b = 3$, while the observed timing ratio is somewhat larger.
Thus, wind alignment provides a substantial computational saving at a modest observed loss in localization performance in this experiment.

\FloatBarrier

\subsection{Baseline Comparisons}
\label{sec:exp_baselines}

This numerical study compares the one-body repeated-cycle controller against two baselines to assess the benefit of belief-driven adaptation and to compare it with a practical myopic information-seeking controller.
The \emph{non-adaptive patrol} controller follows the fixed wind-aligned sinusoidal patrol field from \Cref{eq:patrol_field} and represents a hand-designed crosswind exploration strategy in the spirit of bio-inspired plume tracing \cite{kennedy83zigzagging,li2001crosswindimportance}.
The patrol parameters in \Cref{tab:simulation_scenarios} define the reference field used in the one-body HJB running cost.
For the non-adaptive controller, which follows its reference field for the complete mission without belief feedback, we instead use the domain-scaled parameters $T_p=40\,\mathrm{s}$ and $\rho=0.99$.
They reduce the constant upwind speed to $v_{\max}\sqrt{1-\rho^2}=1.69\,\mathrm{m\,s^{-1}}$ and the maximum crosswind excursion over one period to $\rho v_{\max}T_p/\pi=151\,\mathrm{m}$.
This configuration was selected from the domain geometry to avoid sustained boundary saturation.
The \emph{greedy one-step EIG} controller is a myopic information-theoretic planner that captures the one-step-lookahead philosophy of classical information-driven search \cite{vergassola2007infotaxis,hutchinson2019control}.
Rather than reimplementing those methods, the greedy controller uses the same estimator, plume model, and EIG surrogate as our method.
At each measurement time, it directly solves the static one-step target-selection problem \eqref{eq:greedy_eig}, without constructing or evaluating an HJB value function.
At each planning time $t_k$ it solves
\begin{equation}
\label{eq:greedy_eig}
\bX^\star = \operatorname*{argmax}_{\bX:\,\|\x_i - \x_i(t_k)\|\le v_{\max}\Delta}
\operatorname{EIG}(\bX;\pi_{t_k}^-),
\end{equation}
choosing target positions that maximize the EIG of a single upcoming measurement batch within the reachable ball, then steers each drone straight toward its target at $v_{\max}$ with the same pairwise collision-avoidance correction as the one-body controller \eqref{eq:one_body_collision_correction}.
Problem~\eqref{eq:greedy_eig} is solved by SLSQP with the Euclidean reachable-ball constraints and physical-domain bounds imposed directly, using multiple feasible random restarts; it requires no HJB solve.
The experiment therefore compares the complete control-effort-aware HJB controller with a practical policy based on myopic EIG target selection.

All three controllers are evaluated on the baseline scenario from \Cref{tab:simulation_scenarios}.
For each $\mathrm{SNR}_0\in\{50,100,150,200\}$, we generate $50$ randomized source configurations and run all three methods on the same seeds.
A trial is successful if the final source-location estimate lies within $50\,\mathrm{m}$ of the true source, regardless of early stopping.
% REVIEW B7-completed-cost BEGIN
We report success rate and the cumulative wind-relative control-effort proxy $\mathcal{E}$ from \Cref{sec:setup}.
Final effort includes the stopping rule's effect on run duration; it measures the completed operational policy, not uniform superiority at every common time or effort budget.
% REVIEW B7-completed-cost END
Across the $200$ passive-patrol trials, every trajectory briefly contacts a $y$-boundary and none contacts an $x$-boundary.
The median fraction of mission time with at least one drone on a boundary is $6.6\%$, while the median fraction over all drone-time is $2.2\%$.
Thus, boundary contact is systematic because one initial position lies $100\,\mathrm{m}$ from the upper boundary, but it occupies a small fraction of the simulated trajectories.

\begin{figure}[htbp]
  \centering
  \includegraphics[width=\linewidth]{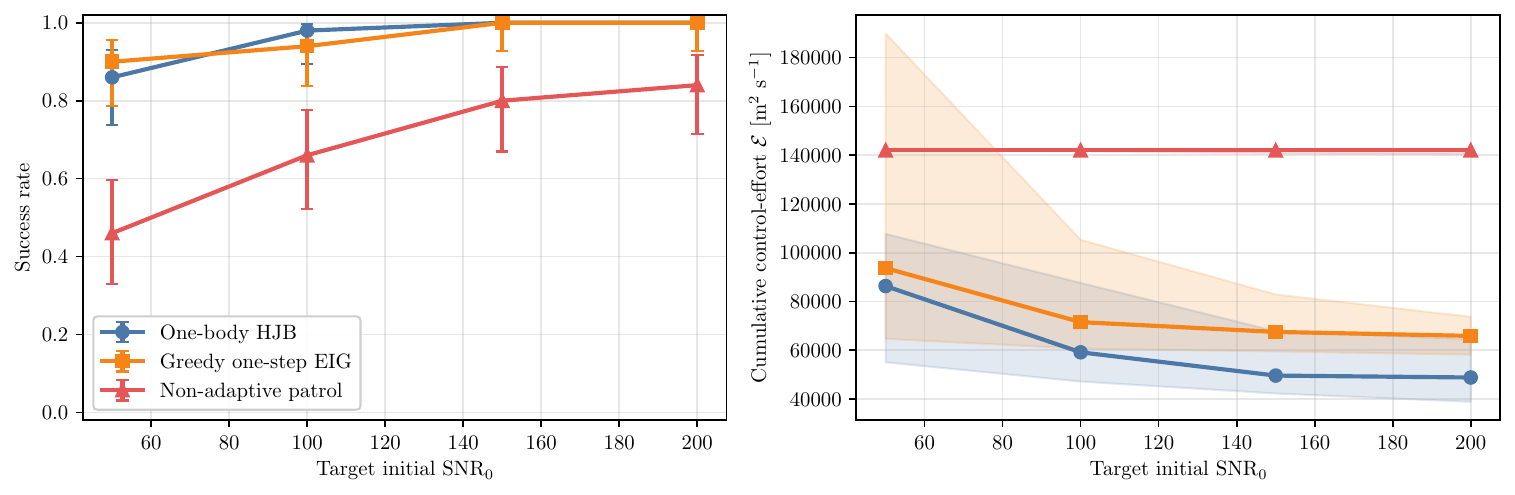}
  \caption{Baseline comparison over $50$ seeds per SNR level in the three-drone baseline scenario.
  Left: success rate versus $\mathrm{SNR}_0$, with Wilson $95\%$ confidence intervals.
  Right: median cumulative control-effort proxy with $10$th--$90$th percentile bands.
  Both adaptive methods achieve higher success rates and lower median control effort than non-adaptive patrol; one-body HJB uses the least control effort among the adaptive controllers.}
  \label{fig:baseline_comparison}
\end{figure}

\Cref{fig:baseline_comparison} summarizes the aggregate outcomes.
Against patrol, the adaptive controllers improve both reported outcomes: at $\mathrm{SNR}_0=50$ the success rate rises from $0.46$ to $0.86$ for one-body HJB and $0.90$ for greedy EIG, while the median control-effort proxy drops from $1.42\times10^5$ to $8.63\times10^4$ and $9.37\times10^4$, respectively.
The patrol success rate reaches $0.84$ at $\mathrm{SNR}_0=200$, whereas both adaptive controllers reach $1.00$ by $\mathrm{SNR}_0=150$.

Against greedy EIG, the one-body HJB controller achieves comparable reliability---the success rates differ by at most four percentage points across the sweep---but consistently uses less control effort.
At $\mathrm{SNR}_0=50$, one-body HJB reduces the median proxy by about $8\%$ relative to greedy EIG; at $\mathrm{SNR}_0=200$ the gap widens to roughly $26\%$ ($4.88\times10^4$ versus $6.59\times10^4$).
Thus, the integrated HJB controller achieves a better empirical reliability--control-effort trade-off than this myopic EIG target-selection policy, especially as the signal level increases.
Because the two controllers differ in more than their planning horizon, this result should not be attributed to multi-step continuation alone.

\begin{figure}[htbp]
  \centering
  \includegraphics[width=.95\linewidth]{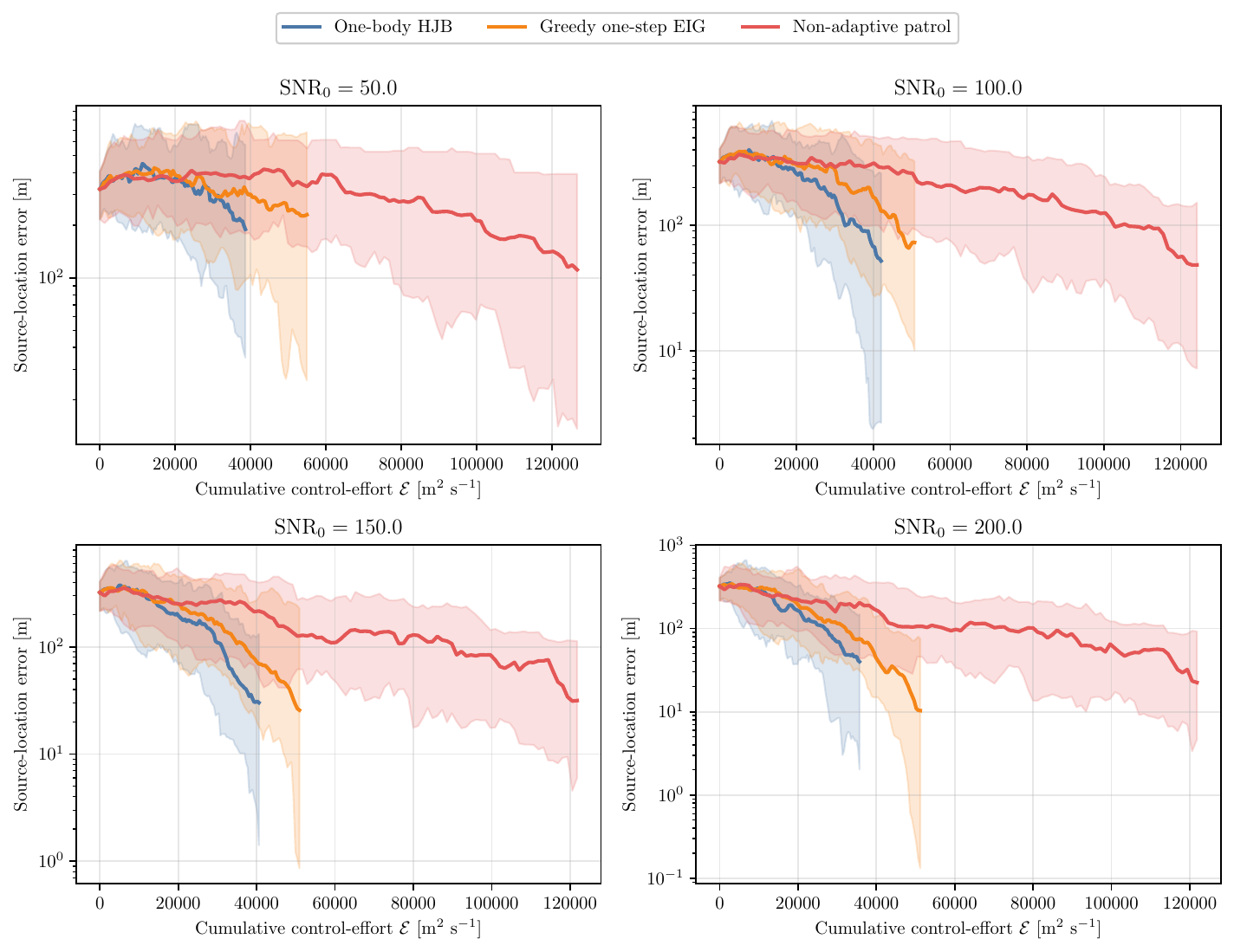}
  \caption{Distance to source ($\|\bs\mu_t-\bs\theta\|$) versus cumulative wind-relative control effort at $\mathrm{SNR}_0=50,100,150,$ and $200$.
% REVIEW B7-curve-support BEGIN
  Solid curves show medians and shaded bands the $10$th--$90$th percentiles over the recorded effort interval shared by all runs within each method; stopped runs are not extended, and support can differ between methods.
% REVIEW B7-curve-support END
  Both adaptive methods bring the source-location estimate closer to the source at lower control effort than patrol, with one-body HJB using less control effort than greedy EIG.}
  \label{fig:baseline_distance_vs_energy}
\end{figure}

\Cref{fig:baseline_distance_vs_energy} examines the accuracy--control-effort trade-off in more detail.
Across the tested signal levels, patrol requires more cumulative control effort to reduce localization error than either adaptive method.
Both adaptive controllers reduce distance to source much earlier along the trajectory, and one-body HJB generally reaches a given error level at lower control effort than greedy EIG.
Together, these results establish that belief-driven replanning outperforms passive crosswind sweeping in this scenario and that the complete HJB controller attains comparable reliability to greedy EIG with lower wind-relative control effort.

\FloatBarrier

\subsection{One-Body vs.\ Pairwise Planning}
\label{sec:exp_pairwise_showcase}

The baseline comparisons in \Cref{sec:exp_baselines} show that one-body HJB replanning outperforms passive patrol and uses less control effort than greedy one-step EIG at comparable localization reliability.
We now compare the complete one-body and pairwise surrogate controllers: does the pairwise construction produce more reliable cooperative localization when both use the same estimator, plume model, no-fly zones, runtime collision safeguard, stopping rule, and discount $r=0.3\,\mathrm{s^{-1}}$?
The pairwise construction changes both the information reward, by retaining two-sensor terms, and the feedback approximation, by considering pair interactions in two-agent HJB subproblems.
All runs use the challenging cooperative scenario from \Cref{tab:simulation_scenarios}, whose geometry is shown in \Cref{fig:pairwise_showcase_setting}.

We run $100$ independent paired trials.
For each seed, the same source location, simulation-noise stream, and MAP candidate sets are used for the one-body and pairwise controllers.
No collision or no-fly-zone violation is observed in these simulated runs.
\Cref{tab:pairwise_showcase_summary} summarizes the aggregate outcomes.

\begin{table}[t]
  \centering
  \small
  \setlength{\tabcolsep}{4pt}
  \caption{Pairwise showcase results over $100$ paired seeds at $r=0.3\,\mathrm{s^{-1}}$.
  Stopped denotes satisfaction of the full diagnostic stopping rule by the time horizon, including the final check.
% REVIEW C2-medians BEGIN
  The success columns report fractions with final error at or below each threshold; error, trace, effort, and runtime are medians, except the explicitly labeled 90th-percentile error.}
% REVIEW C2-medians END
  \label{tab:pairwise_showcase_summary}
  \resizebox{\linewidth}{!}{%
  \begin{tabular}{lcccccccc}
    \toprule
    Method & Stopped & $\le2\,\mathrm{m}$ & $\le5\,\mathrm{m}$ & Error [m] & 90th-percentile error [m] & $\tr(\bs\Sigma)\,[\mathrm{m}^2]$ & $\mathcal E\,[\mathrm{m}^2\mathrm{s}^{-1}]$ & Runtime [s] \\
    \midrule
    One-body & $0.85$ & $0.36$ & $0.71$ % REVIEW C2-fractions
 & $2.92$ & $17.88$ & $12.57$ & $1.54\times 10^4$ & $16.7$ \\
    Pairwise & $0.99$ & $0.69$ & $0.93$ % REVIEW C2-fractions
 & $1.44$ & $4.17$ & $2.70$ & $1.16\times 10^4$ & $263.5$ \\
    \bottomrule
  \end{tabular}%
  }
\end{table}
Pairwise planning produces more reliable localization across the tested trials: the fraction of runs ending within $5\,\mathrm{m}$ of the source rises from $71\%$ to $93\%$.
The lower median and 90th-percentile errors in \Cref{tab:pairwise_showcase_summary} show improvements in both typical outcomes and the upper tail of the error distribution.
Among runs that met the stopping rule, $23$ of $85$ one-body runs and $6$ of $99$ pairwise runs ended more than $5\,\mathrm{m}$ from the source.
% REVIEW B7-runtime BEGIN
The trade-off is computation: the median wall-clock runtime increases from $16.7\,\mathrm{s}$ to $263.5\,\mathrm{s}$ per run, approximately a $16\times$ increase in this experiment.
These whole-run times are not per-cycle latency measurements or evidence of real-time feasibility.
% REVIEW B7-runtime END

\begin{table}[t]
  \centering
  \small
  \setlength{\tabcolsep}{5pt}
  \caption{Paired differences in the pairwise showcase over $100$ paired seeds.
  Differences are one-body minus pairwise, so positive values favor pairwise planning because smaller values are better for all three metrics.
  Intervals are BCa bootstrap $95\%$ confidence intervals.}
  \label{tab:pairwise_showcase_paired_differences}
  \begin{tabular}{lccc}
    \toprule
    Metric & Median difference & Seeds favoring one-body & Seeds favoring pairwise \\
    \midrule
    Error [m] & $1.54\,[0.59,1.99]$ & $33$ & $67$ \\
    $\tr(\bs\Sigma)\,[\mathrm{m}^2]$ & $7.23\,[2.58,12.73]$ & $24$ & $76$ \\
    $\mathcal E\,[\mathrm{m}^2\mathrm{s}^{-1}]$ & $2.95\,[2.56,4.08]\times10^3$ & $13$ & $87$ \\
    \bottomrule
  \end{tabular}
\end{table}

\begin{figure}[htbp]
  \centering
  \includegraphics[width=\linewidth]{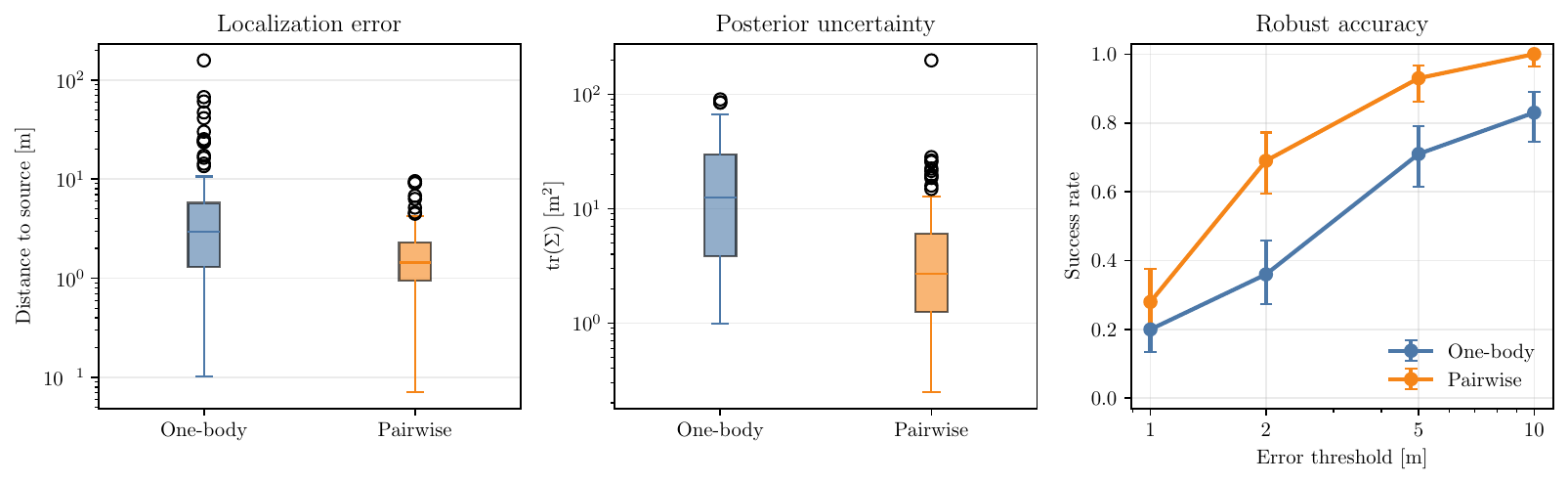}
  \caption{Aggregate pairwise-showcase outcomes over $100$ paired seeds.
  The first two panels show the distributions of final localization error and final posterior trace on log scales.
  The third panel shows success rate as a function of the final-error threshold, with Wilson $95\%$ confidence intervals.}
  \label{fig:pairwise_showcase_summary}
\end{figure}

\Cref{fig:pairwise_showcase_summary} presents the same final-outcome comparison distributionally, including thresholded success rates.
Pairwise planning reaches lower final localization error and covariance trace than one-body planning, while increasing success rate for different error thresholds.
% The effort-normalized trajectory panel in \Cref{fig:pairwise_showcase_energy} does not show uniform pairwise dominance: over the restricted effort range shared by all runs of each method, the one-body median is generally lower.
% The complete-run paired-improvement scatter and \Cref{tab:pairwise_showcase_paired_differences} instead compare each seed against itself: $d_{\mathrm{1b}}-d_{\mathrm{2b}}>0$ and $\mathcal E_{\mathrm{1b}}-\mathcal E_{\mathrm{2b}}>0$ both favor pairwise planning.
% Pairwise planning improves both final localization error and control effort in $58$ of $100$ paired trials, improves localization error only in $9$, improves control effort only in $29$, and is worse in both metrics in $4$.
% The median paired improvement is $1.54\,\mathrm{m}$ in final error and $2.95\times10^3\,\mathrm{m^2s^{-1}}$ in cumulative control effort.

% \begin{figure}[htbp]
%   \centering
%   \begin{subfigure}{0.48\linewidth}
%     \centering
%     \includegraphics[width=\linewidth]{figures/pairwise_showcase_distance_vs_energy.pdf}
%     \caption{Distance to source versus cumulative control effort}
%   \end{subfigure}
%   \begin{subfigure}{0.48\linewidth}
%     \centering
%     \includegraphics[width=\linewidth]{figures/pairwise_showcase_paired_improvement.pdf}
%     \caption{Paired improvements by seed}
%   \end{subfigure}
%   \caption{Control-effort outcomes for the pairwise showcase.
%   In panel (a), each run is interpolated over the effort range shared by all runs of the same method; solid curves show pointwise medians and shaded bands indicate the $10$th--$90$th percentile range.
%   Because this common range excludes later portions of some runs, panel (a) does not compare their final outcomes.
%   In panel (b), positive values on each axis favor pairwise planning for the same seed.}
%   \label{fig:pairwise_showcase_energy}
% \end{figure}

\begin{figure}[htbp]
  \centering
  \begin{subfigure}{0.48\linewidth}
    \centering
    \includegraphics[width=\linewidth]{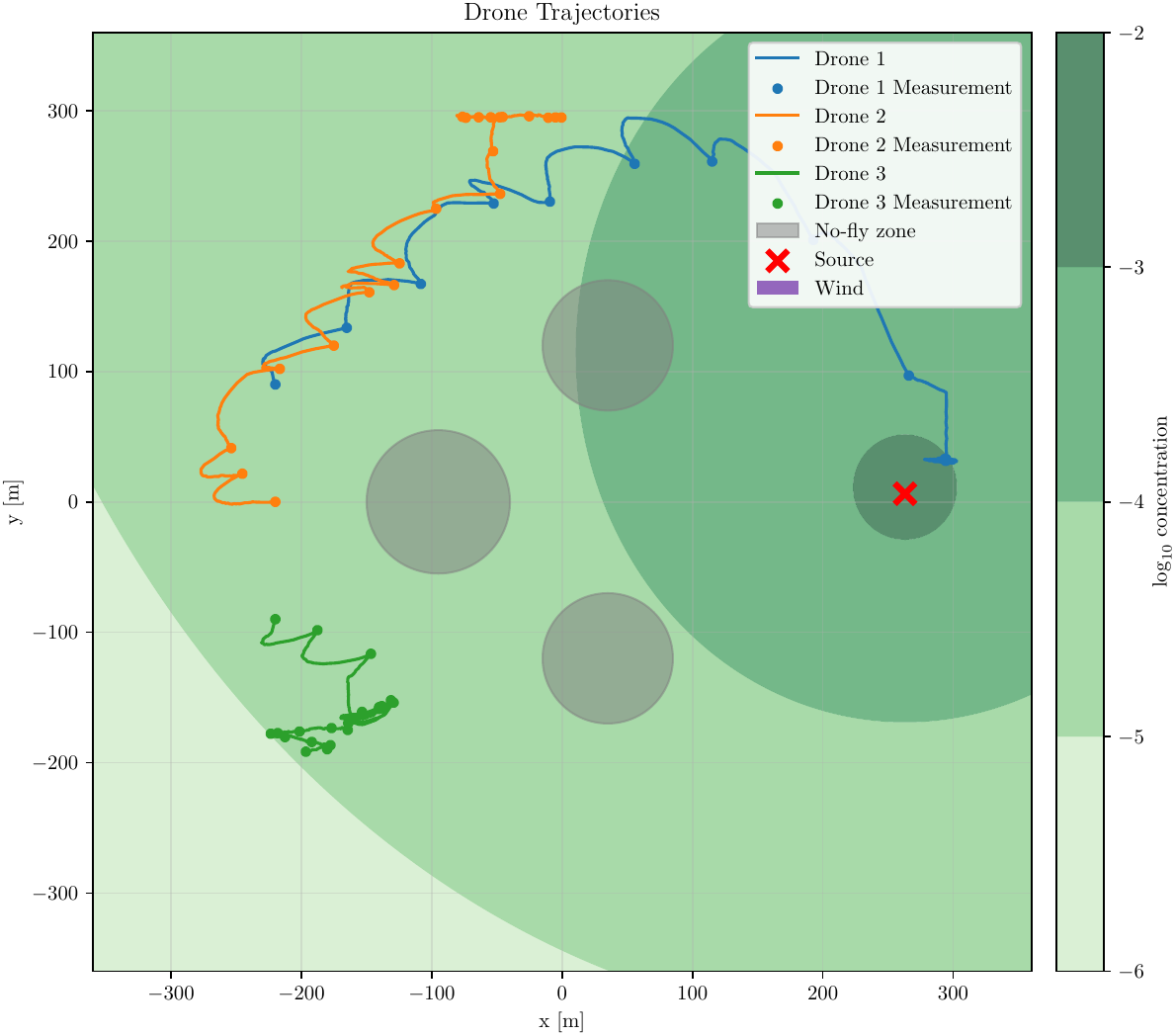}
    \caption{One-body}
  \end{subfigure}
  \begin{subfigure}{0.48\linewidth}
    \centering
    \includegraphics[width=\linewidth]{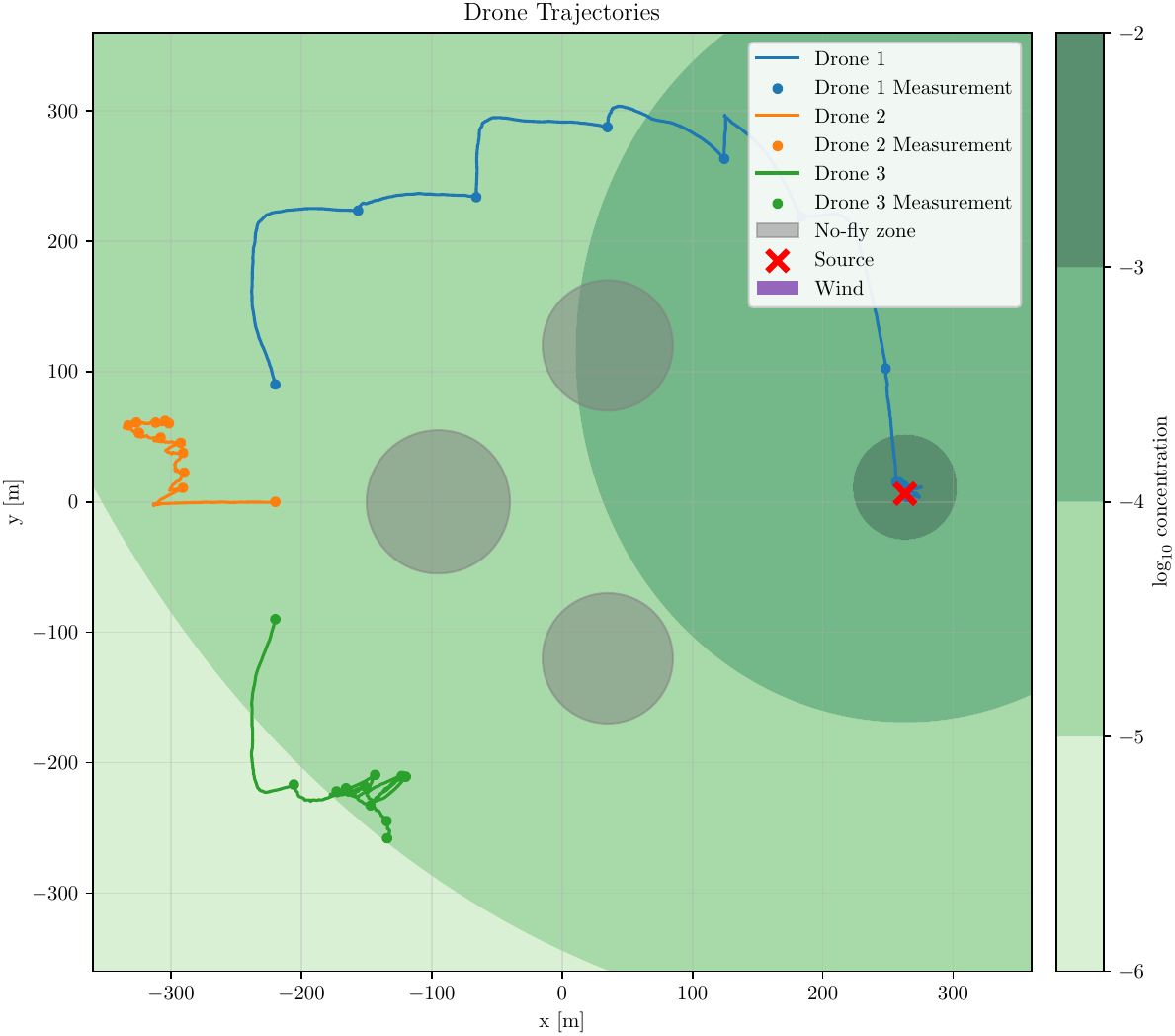}
    \caption{Pairwise}
  \end{subfigure}
  \caption{Illustrative trajectories for one-body (left) and pairwise (right) planning for one paired trial.
  The concentration of pollutant is presented in different shades of green.
  Both panels show the same source and no-fly-zone geometry.
  }
  \label{fig:pairwise_showcase_seed45_paths}
\end{figure}

\begin{figure}[htbp]
  \centering
  \includegraphics[width=.75\linewidth]{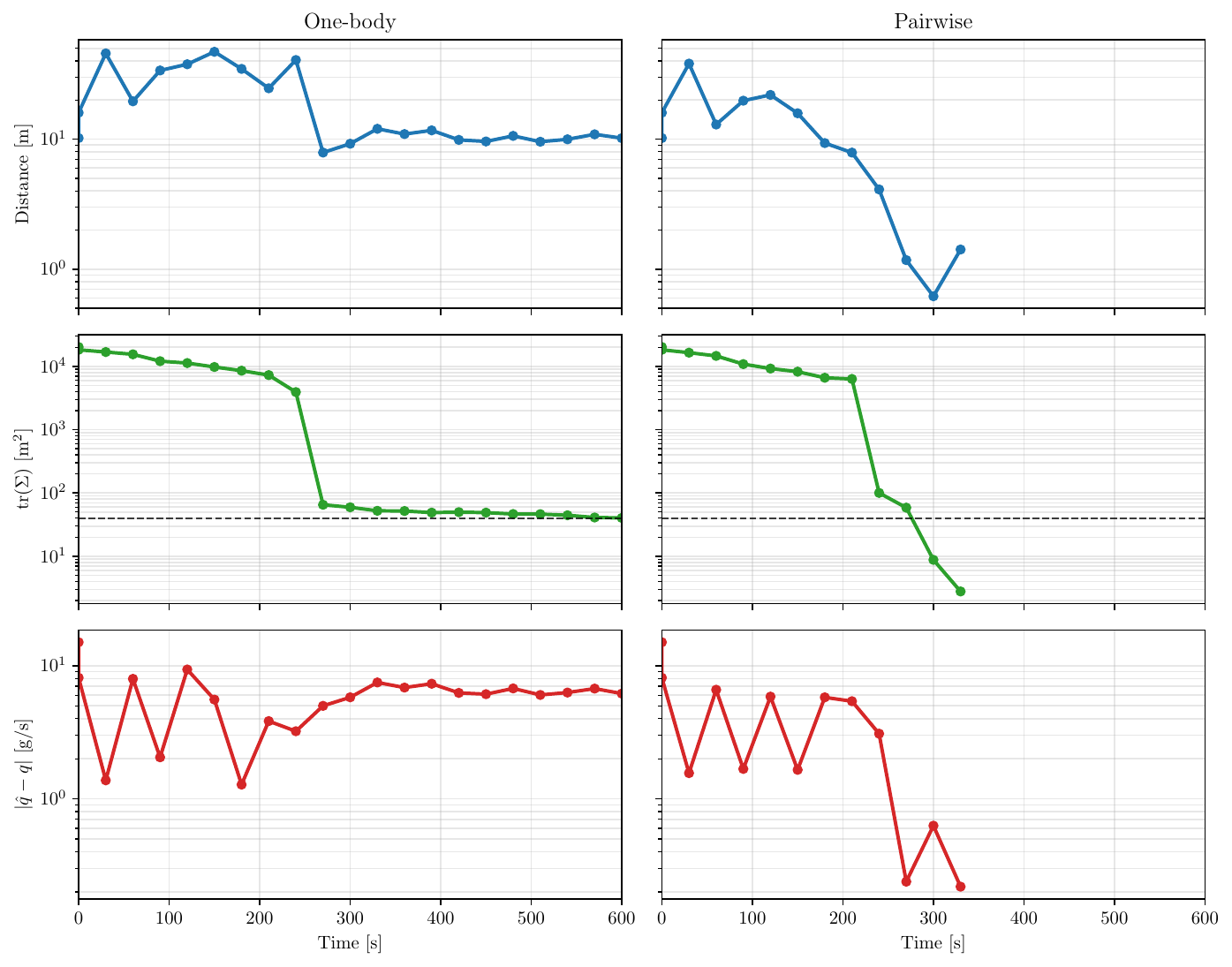}
  \caption{Belief convergence for one-body (left) and pairwise (right) planning for the illustrative paired trial.
  The pairwise controller reduces source-location error, posterior trace, and intensity error earlier, and satisfies the stopping criterion before the horizon.}
  \label{fig:pairwise_showcase_seed45_convergence}
\end{figure}

Finally, \Cref{fig:pairwise_showcase_seed45_paths,fig:pairwise_showcase_seed45_convergence} show an illustrative paired trial.
Here, one-body planning reaches the horizon with final error $10.18\,\mathrm{m}$, $\tr(\bs\Sigma)=40.39\, \mathrm{m}^2$, cumulative control effort $1.85\times10^4\,\mathrm{m^2s^{-1}}$, and travel distance $4248\,\mathrm{m}$.
In contrast, pairwise planning satisfies the diagnostic stopping criterion at $t=330\,\mathrm{s}$, with final error $1.42\,\mathrm{m}$, $\tr(\bs\Sigma)=2.79\, \mathrm{m}^2$, cumulative control effort $1.12\times10^4\,\mathrm{m^2s^{-1}}$, and travel distance $2410\,\mathrm{m}$.
The path plots suggest more complementary sensing paths around the obstacle field under pairwise feedback, followed by rapid movement of the source-location estimate toward the source.
This single trial illustrates a possible coordination mechanism; the aggregate comparisons establish the difference in localization outcomes.
The convergence plots show that this is accompanied by simultaneous reduction of source-location error, posterior trace, and intensity error.

\FloatBarrier

\subsection{Discount Sensitivity}
\label{sec:exp_discount_sensitivity}

The previous experiment fixes the discount rate at $r=0.3\,\mathrm{s^{-1}}$.
We now ask how sensitive the planning is to this choice.
In the Bellman fixed point \eqref{eq:bellman}, $r$ discounts both the next information reward and the continuation value relative to costs incurred earlier in the cycle.
A larger $r$ therefore makes the policy more myopic: distant benefits of moving into informative configurations are discounted more heavily, and the controller should be less willing to spend control effort on exploratory positioning.
The corresponding dimensionless per-cycle factor is $\beta=e^{-r\Delta}$, which multiplies the complete measurement-time jump term $V(0,S_\Delta^+)-\operatorname{EIG}(\bX_\Delta;\pi_{t_k}^-)$ in \eqref{eq:bellman}.

We use the challenging cooperative scenario from \Cref{tab:simulation_scenarios}, with $\Delta=30\,\mathrm{s}$, and run the same $100$ paired seeds for one-body and pairwise planning at
\begin{equation}
r\in\{0.1,0.2,0.3,0.5,0.7\}\,\mathrm{s^{-1}}.
\end{equation}
The associated per-cycle factors are listed in \Cref{tab:discount_factors}.
\begin{table}[h]
  \centering
  \caption{Discount rates and corresponding per-cycle factors for $\Delta=30\,\mathrm{s}$.}
  \label{tab:discount_factors}
  \begin{tabular}{cc}
    \toprule
    Discount rate $r\,[\mathrm{s^{-1}}]$ & Per-cycle factor $\beta=e^{-30r}$ \\
    \midrule
    $0.1$ & $4.98\times10^{-2}$ \\
    $0.2$ & $2.48\times10^{-3}$ \\
    $0.3$ & $1.23\times10^{-4}$ \\
    $0.5$ & $3.06\times10^{-7}$ \\
    $0.7$ & $7.58\times10^{-10}$ \\
    \bottomrule
  \end{tabular}
\end{table}

Relative to $r=0.3\,\mathrm{s^{-1}}$, the rates $r=0.1$ and $0.2\,\mathrm{s^{-1}}$ increase $\beta$ by factors of approximately $400$ and $20$, respectively, while preserving the physical and numerical setting.
For each seed and rate, the methods use the same source location and noise realizations.
No collision or no-fly-zone violation is observed in these simulated runs.
All accepted fixed-point solutions satisfy the normalized Bellman-residual criterion in \eqref{eq:normalized_bellman_residual}.
At $r=0.1\,\mathrm{s^{-1}}$, nine pairwise trials require at least one successful retry with $\Delta t_{\mathrm{HJB}}/2$.
\Cref{tab:discount_sweep_summary} summarizes the results.

\begin{table}[t]
  \centering
  \small
  \setlength{\tabcolsep}{4pt}
  \caption{Discount-sensitivity results in the challenging cooperative scenario.
  Stopped denotes satisfaction of the diagnostic stopping rule by the time horizon.
  Error, $\tr(\bs\Sigma)$, control effort, travel distance, and runtime are medians.
  Every row contains $100$ trials.}
  \label{tab:discount_sweep_summary}
  \begin{tabular}{c l c c c c c c c}
    \toprule
    $r\,[\mathrm{s^{-1}}]$ & Method & Stopped & $\le5\,\mathrm{m}$ & Error [m] & $\tr(\bs\Sigma)\,[\mathrm{m}^2]$ & $\mathcal E\,[\mathrm{m}^2\mathrm{s}^{-1}]$ & Distance [m] & Runtime [s] \\
    \midrule
    \multirow{2}{*}{0.1}
      & One-body & $0.98$ & $0.93$ & $1.10$ & $2.14$ & $1.99\times10^4$ & $2825$ & $8.5$ \\
      & Pairwise & $1.00$ & $0.94$ & $0.94$ & $1.68$ & $2.16\times10^4$ & $2920$ & $178.1$ \\
    \midrule
    \multirow{2}{*}{0.2}
      & One-body & $0.97$ & $0.76$ & $1.44$ & $4.08$ & $1.68\times10^4$ & $2861$ & $12.6$ \\
      & Pairwise & $1.00$ & $0.93$ & $1.30$ & $2.58$ & $1.15\times10^4$ & $2251$ & $214.8$ \\
    \midrule
    \multirow{2}{*}{0.3}
      & One-body & $0.85$ & $0.71$ & $2.92$ & $12.57$ & $1.54\times10^4$ & $3186$ & $16.7$ \\
      & Pairwise & $0.99$ & $0.93$ & $1.44$ & $2.70$ & $1.16\times10^4$ & $2431$ & $263.5$ \\
    \midrule
    \multirow{2}{*}{0.5}
      & One-body & $0.73$ & $0.72$ & $2.74$ & $9.09$ & $1.36\times10^4$ & $3513$ & $25.0$ \\
      & Pairwise & $0.96$ & $0.93$ & $1.24$ & $2.95$ & $1.20\times10^4$ & $3029$ & $376.0$ \\
    \midrule
    \multirow{2}{*}{0.7}
      & One-body & $0.13$ & $0.32$ & $20.41$ & $1.42\times10^3$ & $1.20\times10^4$ & $3551$ & $27.4$ \\
      & Pairwise & $0.34$ & $0.56$ & $3.99$ & $10.84$ & $1.18\times10^4$ & $3488$ & $473.1$ \\
    \bottomrule
  \end{tabular}
\end{table}

\begin{figure}[htbp]
  \centering
  \includegraphics[width=\linewidth]{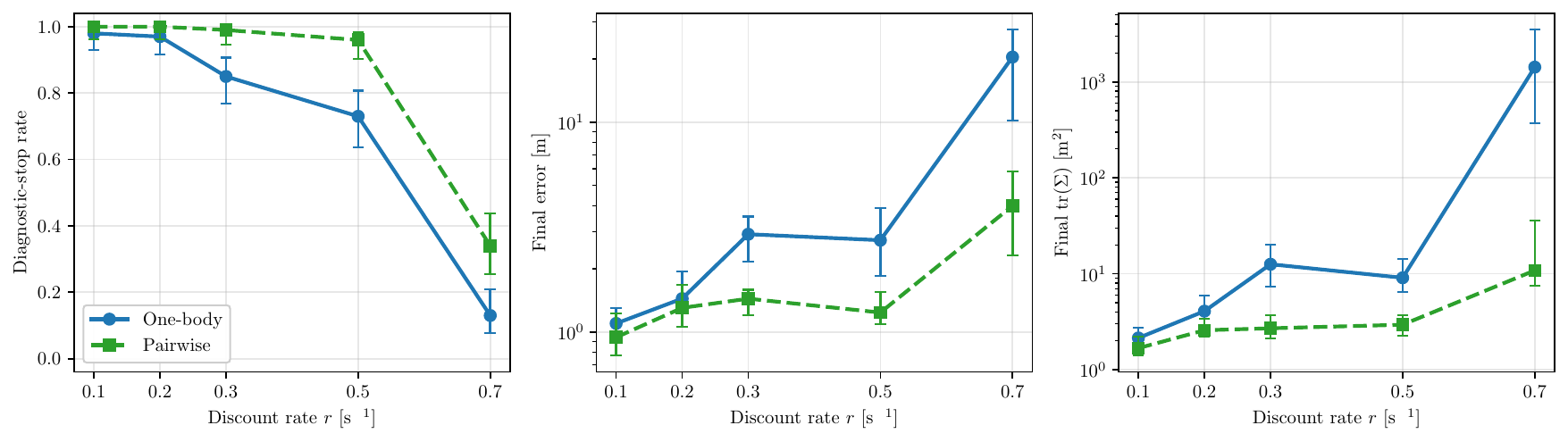}
  \caption{Localization outcomes as a function of discount rate $r$; the corresponding per-cycle factors are listed in \Cref{tab:discount_factors}.
  The left panel shows diagnostic-stop rates with Wilson $95\%$ intervals.
  The middle and right panels show medians with BCa bootstrap $95\%$ intervals.
  Every point uses $100$ trials.}
  \label{fig:discount_sweep_outcomes}
\end{figure}

\begin{figure}[htbp]
  \centering
  \includegraphics[width=\linewidth]{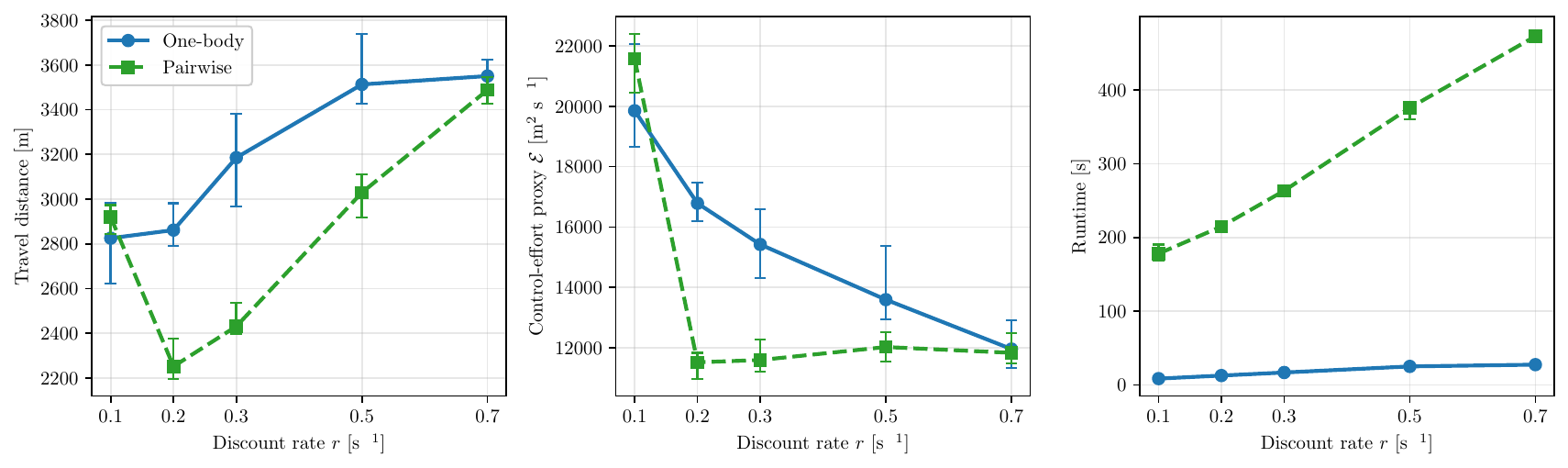}
  \caption{Travel distance, cumulative wind-relative control-effort proxy, and wall-clock runtime as functions of $r$.
  Markers show medians and error bars show BCa bootstrap $95\%$ intervals over $100$ trials per point.}
  \label{fig:discount_sweep_costs}
\end{figure}

Pairwise planning maintains reliable localization over a wider range of discount rates than one-body planning.
Pairwise localization remains reliable through $r=0.5\,\mathrm{s^{-1}}$, whereas the one-body outcomes vary more strongly across the intermediate rates.
At $r=0.7\,\mathrm{s^{-1}}$, both methods deteriorate, with the failure especially pronounced for one-body planning, whose median final error rises to $20.41\,\mathrm{m}$.
This deterioration is consistent with heavy discounting reducing the incentive to move toward informative measurement locations.
Because the next information reward and continuation value are discounted together, the sweep does not identify which component drives the change.

For pairwise planning, reducing $r$ from $0.3$ to $0.2\,\mathrm{s^{-1}}$ reduces control effort and travel distance, while the paired changes in error and posterior trace are not statistically resolved.
Reducing the rate further to $0.1\,\mathrm{s^{-1}}$ decreases the pairwise median error by $28\%$, but increases median control effort by $87\%$ and median travel distance by $30\%$.
The same qualitative accuracy--effort trade-off appears for one-body planning.

Pairwise planning is substantially more expensive computationally throughout the sweep, so its wider range of reliable discount rates comes at a clear runtime cost.
% REVIEW B2-discount BEGIN
Although the largest tested per-cycle factor is only $0.0498$, the observed changes establish sensitivity to discounting of the complete measurement-time jump, which scales immediate EIG and continuation together.
At these discounts, a small Bellman residual does not by itself demonstrate a substantial contribution from continuation.
% REVIEW B2-discount END

\FloatBarrier

\subsection{Scaling with Team Size}
\label{sec:exp_scaling_M}

Both the one-body and pairwise planning surrogates keep the per-cycle HJB subproblem at a fixed dimension independent of team size $M$: one-body planning always solves a 5-dimensional subproblem (two position coordinates plus three belief coordinates), while pairwise planning solves a 7-dimensional one.
One-body and pairwise table queries cost $O(M)$ and $O(M^2)$, respectively, while the all-pairs collision correction makes both complete commands $O(M^2)$ in the present implementation.
We use the pairwise surrogate because it produced more reliable localization in the preceding cooperative comparison.
We examine how localization accuracy and control effort per drone change with team size, using the challenging cooperative scenario from \Cref{tab:simulation_scenarios} with $M\in\{2,3,4\}$.

\begin{figure}[t]
  \centering
  \includegraphics[width=\linewidth]{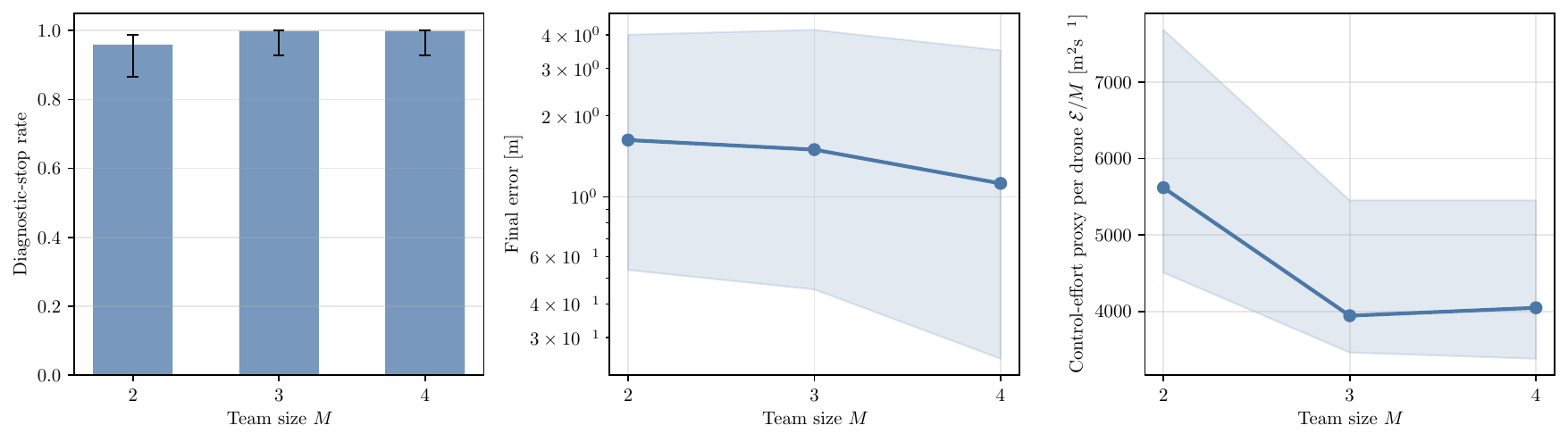}
  \caption{Team-size scaling over $50$ seeds with the pairwise surrogate, $M\in\{2,3,4\}$.
  Left: diagnostic-stop rate (fraction of trials that meet the stopping rule by the time horizon), with Wilson $95\%$ confidence intervals.
  Center: distribution of final localization error (median line, 10th--90th percentile band).
  Right: distribution of the control-effort proxy per drone $\mathcal{E}/M$ (same bands).}
  \label{fig:scaling_M}
\end{figure}

For each $M$, we use the same source box, prior, plume model, no-fly zones, collision model, and stopping criterion ($r=0.3\,\mathrm{s^{-1}}$, $50$ seeds).
The starting positions are evenly spaced along $x=-220\,\mathrm{m}$: $(-220,\pm90)$ for $M=2$; $(-220,\,90)$, $(-220,\,0)$, $(-220,\,-90)$ for $M=3$; and $(-220,\,90)$, $(-220,\,30)$, $(-220,\,-30)$, $(-220,\,-90)$ for $M=4$.
We report $\mathcal{E}/M$ to compare control effort per drone across team sizes.
\Cref{fig:scaling_M} summarizes the resulting stop rate, final localization error, and per-drone control effort as $M$ increases.

Localization reliability is already near saturation at $M=2$ and remains high for the larger teams.
The $5\,\mathrm{m}$ success rate is not monotone in $M$, however, so the experiment does not support a uniform improvement in every thresholded outcome.
The clearer trend is in median final error, which decreases from $1.63\,\mathrm{m}$ at $M=2$ to $1.12\,\mathrm{m}$ at $M=4$.

The control effort per drone falls by roughly $30\%$ from $M=2$ to $M=3$ and then changes little at $M=4$.
Because total team effort increases with team size, this pattern indicates a redistribution of effort across more agents rather than a reduction in aggregate mission cost.
Thus, the experiment shows that the fixed 7-D pairwise HJB subproblem can be reused as the team grows from two to four drones while maintaining reliable localization and avoiding growth in per-agent effort.

\FloatBarrier

\section{Conclusion}
\label{sec:conclusion}

We presented a framework for multi-UAV gas source localization that combines continuous HJB feedback with discrete Bayesian updates of source location and intensity.
A repeated-cycle Bellman problem balances motion cost against expected information gain, while reduced belief coordinates and reusable one-body or pairwise subproblems limit the dimension of the HJB solve.

The experiments show that the choice of planning surrogate trades computational cost against localization performance.
One-body HJB control uses less control effort than patrol and greedy EIG in the baseline scenario, with success rates close to those of greedy EIG.
Pairwise planning improves localization reliability in the challenging cooperative scenario, but requires substantially more computation.
Wind alignment reduces planning time at a modest observed loss in localization performance in the one-body comparison, while the team-size study demonstrates reuse of the same pairwise HJB subproblem for two to four drones.
Together, these results support using reduced HJB feedback to balance sensing performance, motion cost, and available computation.
The discount sweep exposes a further accuracy--effort trade-off: weaker discounting can improve localization at increased motion cost.

The main limitations arise from the approximations used for tractability and from the scope of the simulations.
The block-diagonal Fisher-information model and the one-body and pairwise feedback constructions are not error-controlled approximations of the full joint belief-state problem.
Moreover, the experiments use a steady plume surrogate, a proxy for propulsion energy, known wind, homogeneous sensors, centralized estimation, and perfect communication.
Future work should focus on quantifying the planner approximations, testing the controller with time-varying dispersion and imperfect communication in high-fidelity or field experiments, and evaluating larger teams with a decomposition of HJB solution, feedback evaluation, and data-assimilation costs.

\section*{Acknowledgments and Funding}

This work was partly performed as part of the Helmholtz School for Data Science in Life, Earth and Energy (HDS-LEE) and received funding from the Helmholtz Association of German Research Centres.
André Carlon carried out most of his contribution while a researcher at the Chair of Mathematics for Uncertainty Quantification, RWTH Aachen University, with funding from the Alexander von Humboldt Foundation.
Raúl Tempone acknowledges support from an Alexander von Humboldt Professorship, KAUST baseline funding, and KAUST award CRG2025-1862, ``Efficient Hierarchical Approximations for Partially Observed Stochastic Optimal Control with Applications.''

\appendix

\section{Expected Information Gain for the Linearized-Gaussian Model}
\label{app:eig}
This appendix derives the linearized-Gaussian EIG used in \Cref{prop:eig}.
We first derive the location component of the block-diagonal planner surrogate, then give the one-dimensional intensity analogue, and finally combine the two components under approximations (B1)--(B2).

\paragraph{Linearized location model.}
Fix one planning cycle and consider the location subproblem induced by (B2), holding the log-intensity fixed at its current nominal value $\zeta=\hat\zeta_t$ and linearizing at $\bs\theta=\bs\mu_t$.
For a candidate team configuration $\bX=(\x_1,\dots,\x_M)$, the next measurement batch $Y\in\R^M$ is approximated by
\begin{equation}
\label{eq:app_linear_obs}
Y \approx h_0(\bX) + H(\bX)(\bs\theta-\bs\mu_t) + \varepsilon,
\qquad
\varepsilon\sim\mathcal N(0,\sigma_e^2 \bs I_M),
\end{equation}
where
\begin{equation}
h_0(\bX)=
\begin{bmatrix}
C(\x_1,\bs\mu_t,\hat\zeta_t)\\[-1mm]
\vdots\\[-1mm]
C(\x_M,\bs\mu_t,\hat\zeta_t)
\end{bmatrix},
\qquad
H(\bX)=
\begin{bmatrix}
\nabla_{\bs\theta}C(\x_1,\bs\mu_t,\hat\zeta_t)^\top\\[-1mm]
\vdots\\[-1mm]
\nabla_{\bs\theta}C(\x_M,\bs\mu_t,\hat\zeta_t)^\top
\end{bmatrix}\in\R^{M\times 2}.
\end{equation}
The cycle-start marginal prior for the location parameter is
\begin{equation}
\pi_{\bs\theta,t}^-
=\mathcal N(\bs\mu_t,\bs\Sigma_t).
\end{equation}

\paragraph{Posterior moments under the Gaussian surrogate.}
Marginalizing $\bs\theta$ under the prior gives the predictive distribution
\begin{equation}
Y\mid\bX,\pi_{\bs\theta,t}^-
\sim\mathcal N\!\bigl(h_0(\bX),\bs S_Y\bigr),
\qquad
\bs S_Y:=H(\bX)\bs\Sigma_tH(\bX)^\top+\sigma_e^2\bs I_M.
\end{equation}
For the linear-Gaussian model \eqref{eq:app_linear_obs}, the location posterior is Gaussian with moments available in closed-form~\cite{kaipio2005statistical}.
The location posterior is
\begin{equation}
\pi_{\bs\theta,t}^+(\cdot\mid Y,\bX)
=\mathcal N(\bs\mu_+,\bs\Sigma_+),
\end{equation}
with mean
\begin{equation}
\label{eq:app_post_mean}
\bs\mu_+
=\bs\mu_t+\bs K\bigl(Y-h_0(\bX)\bigr),
\qquad
\bs K:=\bs\Sigma_tH(\bX)^\top\bs S_Y^{-1},
\end{equation}
and covariance
\begin{equation}
\label{eq:app_post_cov}
\begin{aligned}
\bs\Sigma_+
&=
\left(
\bs\Sigma_t^{-1} + H(\bX)^\top (\sigma_e^2 \bs I_M)^{-1} H(\bX)
\right)^{-1}
\\
&=
\left(
\bs\Sigma_t^{-1}+\frac{1}{\sigma_e^2}H(\bX)^\top H(\bX)
\right)^{-1}
\\
&=
\bs\Sigma_t-
\bs\Sigma_tH(\bX)^\top\bs S_Y^{-1}H(\bX)\bs\Sigma_t.
\end{aligned}
\end{equation}
The last equality is the Woodbury matrix identity.
Since
\begin{equation}
\frac{1}{\sigma_e^2}H(\bX)^\top H(\bX)
=
\frac{1}{\sigma_e^2}\sum_{i=1}^M
\nabla_{\bs\theta}C(\x_i,\bs\mu_t,\hat\zeta_t)\,
\nabla_{\bs\theta}C(\x_i,\bs\mu_t,\hat\zeta_t)^\top,
\end{equation}
we recover exactly the Fisher matrix $\bs F_{\bs\theta}(\bX)$ from \eqref{eq:fisher_theta}, and therefore
\begin{equation}
\label{eq:app_post_cov_fisher}
\bs\Sigma_+ = \bigl(\bs\Sigma_t^{-1}+\bs F_{\bs\theta}(\bX)\bigr)^{-1}.
\end{equation}

\paragraph{Gaussian KL divergence.}
For Gaussian measures $\mathcal N(\bs m_1,\bs\Sigma_1)$ and $\mathcal N(\bs m_2,\bs\Sigma_2)$ in dimension $d$,
\begin{equation}
\label{eq:app_gauss_kl}
\KL{\mathcal N(\bs m_1,\bs\Sigma_1)}{\mathcal N(\bs m_2,\bs\Sigma_2)}
=
\tfrac12\Bigl[
\log\det(\bs\Sigma_2\bs\Sigma_1^{-1})
-d
+\tr(\bs\Sigma_2^{-1}\bs\Sigma_1)
+(\bs m_1-\bs m_2)^\top \bs\Sigma_2^{-1}(\bs m_1-\bs m_2)
\Bigr].
\end{equation}
Applying \eqref{eq:app_gauss_kl} with $(\bs m_1,\bs\Sigma_1)=(\bs\mu_+,\bs\Sigma_+)$ and $(\bs m_2,\bs\Sigma_2)=(\bs\mu_t,\bs\Sigma_t)$ yields
\begin{equation}
\label{eq:app_kl_post_prior}
\KL{\pi_{\bs\theta,t}^+(\cdot\mid Y,\bX)}{\pi_{\bs\theta,t}^-}
=
\tfrac12\Bigl[
\log\det(\bs\Sigma_t\bs\Sigma_+^{-1})
-2
+\tr(\bs\Sigma_t^{-1}\bs\Sigma_+)
+(\bs\mu_+-\bs\mu_t)^\top \bs\Sigma_t^{-1}(\bs\mu_+-\bs\mu_t)
\Bigr].
\end{equation}

\paragraph{Expectation over the incoming measurement batch.}
The location EIG is the expectation of \eqref{eq:app_kl_post_prior} with respect to the prior predictive distribution of $Y$ conditional on $\bX$ and $\pi_{\bs\theta,t}^-$.
For fixed $\bX$, equation \eqref{eq:app_post_cov} shows that $\bs\Sigma_+$ does not depend on the realized value of $Y$; only $\bs\mu_+$ is random before the batch is observed.
From the predictive distribution and \eqref{eq:app_post_mean},
\begin{equation}
\E_{Y\mid\bX,\pi_{\bs\theta,t}^-}[\bs\mu_+]
=\bs\mu_t
\end{equation}
and
\begin{align}
&\E_{Y\mid\bX,\pi_{\bs\theta,t}^-}\!\left[
(\bs\mu_+-\bs\mu_t)(\bs\mu_+-\bs\mu_t)^\top
\right]\\
&\qquad
=\bs K\bs S_Y\bs K^\top
=\bs\Sigma_tH(\bX)^\top\bs S_Y^{-1}H(\bX)\bs\Sigma_t
=\bs\Sigma_t-\bs\Sigma_+,
\end{align}
where the final equality follows from the Woodbury form in \eqref{eq:app_post_cov}.
Equivalently,
\begin{equation}
\bs\Sigma_t
=
\E_{Y\mid\bX,\pi_{\bs\theta,t}^-}\!\left[
(\bs\mu_+-\bs\mu_t)(\bs\mu_+-\bs\mu_t)^\top
\right]
+\bs\Sigma_+.
\end{equation}
This is the law of total covariance written explicitly for the linear-Gaussian surrogate: the prior covariance is the sum of the covariance remaining after the measurement and the covariance, over prospective measurements, of the posterior mean.
Taking the trace after multiplication by $\bs\Sigma_t^{-1}$ yields
\begin{equation}
\label{eq:app_mean_identity}
\E_{Y\mid\bX,\pi_{\bs\theta,t}^-}\!\left[
(\bs\mu_+-\bs\mu_t)^\top \bs\Sigma_t^{-1}(\bs\mu_+-\bs\mu_t)
\right]
=
\tr\!\left(\bs\Sigma_t^{-1}(\bs\Sigma_t-\bs\Sigma_+)\right)
=
2-\tr(\bs\Sigma_t^{-1}\bs\Sigma_+).
\end{equation}
Substituting \eqref{eq:app_mean_identity} into the expectation of \eqref{eq:app_kl_post_prior} cancels the last three terms and gives
\begin{equation}
\label{eq:app_eig_logdet_ratio}
\operatorname{EIG}_{\bs\theta}(\bX;\pi_{\bs\theta,t}^-)
=
\E_{Y\mid\bX,\pi_{\bs\theta,t}^-}\!\left[
\KL{\pi_{\bs\theta,t}^+(\cdot\mid Y,\bX)}{\pi_{\bs\theta,t}^-}
\right]
=
\tfrac12\log\det(\bs\Sigma_t\bs\Sigma_+^{-1}).
\end{equation}

\paragraph{Expression through the Fisher matrix.}
Using \eqref{eq:app_post_cov_fisher},
\begin{equation}
\bs\Sigma_t\bs\Sigma_+^{-1}
=
\bs\Sigma_t\bigl(\bs\Sigma_t^{-1}+\bs F_{\bs\theta}(\bX)\bigr)
=
\bs I_2+\bs\Sigma_t\bs F_{\bs\theta}(\bX).
\end{equation}
By Sylvester's determinant identity,
\begin{equation}
\det\!\bigl(\bs I_2+\bs\Sigma_t\bs F_{\bs\theta}(\bX)\bigr)
=
\det\!\bigl(\bs I_2+\bs\Sigma_t^{1/2}\bs F_{\bs\theta}(\bX)\bs\Sigma_t^{1/2}\bigr).
\end{equation}
Therefore, using the manifestly symmetric form,
\begin{equation}
\label{eq:app_eig_theta_final}
\operatorname{EIG}_{\bs\theta}(\bX;\pi_{\bs\theta,t}^-)
=
\tfrac12\log\det\!\bigl(
\bs I_2+\bs\Sigma_t^{1/2}\bs F_{\bs\theta}(\bX)\bs\Sigma_t^{1/2}
\bigr),
\end{equation}
which is the location term in \eqref{eq:eig_formula}.

\paragraph{One-dimensional intensity term.}
The intensity contribution follows by the same argument for the marginal prior
$\pi_{\zeta,t}^-=\mathcal N(\hat\zeta_t,\sigma_{\zeta,t}^2)$, holding $\bs\theta=\bs\mu_t$.
Its scalar Fisher information is
\begin{equation}
F_\zeta(\bX)
=
\frac{1}{\sigma_e^2}\sum_{i=1}^M
\left(\frac{\partial C}{\partial \zeta}(\x_i,\bs\mu_t,\hat\zeta_t)\right)^2.
\end{equation}
The corresponding posterior variance is
\begin{equation}
\sigma_{\zeta,+}^2=\bigl(\sigma_{\zeta,t}^{-2}+F_\zeta(\bX)\bigr)^{-1},
\end{equation}
and the EIG is
\begin{equation}
\label{eq:app_eig_zeta_final}
\operatorname{EIG}_{\zeta}(\bX;\pi_{\zeta,t}^-)
=
\tfrac12\log\!\bigl(1+\sigma_{\zeta,t}^2F_\zeta(\bX)\bigr).
\end{equation}

\paragraph{Block-diagonal joint surrogate.}
Under (B1)--(B2), define the surrogate prior covariance and Fisher information by
\begin{equation}
\bs\Xi=\diag(\bs\Sigma_t,\sigma_{\zeta,t}^2),
\qquad
\bs F_{\mathrm{BD}}(\bX)=\diag\!\bigl(\bs F_{\bs\theta}(\bX),F_\zeta(\bX)\bigr).
\end{equation}
Applying the Gaussian EIG identity to this deliberately block-diagonal surrogate gives
\begin{align}
\operatorname{EIG}_{\mathrm{BD}}(\bX;\pi_t^-)
&=
\tfrac12\log\det\!\left(
\bs I_3+\bs\Xi^{1/2}\bs F_{\mathrm{BD}}(\bX)\bs\Xi^{1/2}
\right)\\
&=
\tfrac12\log\det\!\left(
\bs I_2+\bs\Sigma_t^{1/2}\bs F_{\bs\theta}(\bX)\bs\Sigma_t^{1/2}
\right)
+\tfrac12\log\!\left(1+\sigma_{\zeta,t}^2F_\zeta(\bX)\right).
\end{align}
Using Sylvester's identity in the location term yields exactly \eqref{eq:eig_formula}.

\bibliographystyle{plain}
\bibliography{references}

\end{document}